\ifx\fltex\undefined
\documentclass[11pt,a4paper,english]{article}
\usepackage{amsmath,amssymb}
\usepackage{babel}
\usepackage{eucal}
\fi 

\DeclareMathVersion{scriptstyle}

\usepackage{isolatin1,theorem}
\theoremstyle{change}
{\theorembodyfont{\slshape} 
  \newtheorem{thm}{Theorem.}[section]
  \newtheorem{lemma}[thm]{Lemma.}
  \newtheorem{prop}[thm]{Proposition.}
  
}
{\theorembodyfont{\rmfamily} 
  \newtheorem{definition}[thm]{Definition.}
  \newtheorem{remark}[thm]{Remark.}

  \newtheorem{notation}[thm]{Notation.}
}

\numberwithin{equation}{section} 

\newcommand{\proof}[1][Proof. ]{{\it#1}}

\def\endproof{{\nobreak\qquad{\scriptstyle \blacksquare}}}

\def\C{{\mathbb C}}
\def\B{{\mathbb B}}

\def\K{{\mathbb K}}

\def\N{{\mathbb N}}
\def\R{{\mathbb R}}

\def\CA{{\mathcal A}}

\def\CD{{\mathcal D}}

\def\CF{{\mathcal F}}
\def\CG{{\mathcal G}}
\def\CH{{\mathcal H}}
\def\CI{{\mathcal I}}

\def\CK{{\mathcal K}}
\def\CL{{\mathcal L}}
\def\CM{{\mathcal M}}

\def\CP{{\mathcal P}}

\def\unit{{\bf 1}}
\def\i{{\rm i}}
\def\e{{\rm e}}
\def\d{{\rm d}}
\def\eps{\varepsilon}
\def\<{{\langle}}
\def\>{{\rangle}}

\def\tr{{\text{\rm tr}}}

\def\im{{\text{\rm Im}}}
\def\re{{\text{\rm Re}}}
\def\ker{{\text{\rm ker}}}
\def\supp{{\text{\rm supp}}}

\def\range{{\text{\rm range}}}

\def\cc{^*}

\begin{document}

\title{{\LARGE{\bf Brown measures of sets of commuting \\operators in a type II$_1$ factor.}}}
\author{Hanne~Schultz \footnote{ As a student of the
      PhD-school OP-ALG-TOP-GEO the author is partially supported by the Danish Research Training Council.} \footnote{Partially
      supported by The Danish National Research Foundation.} }
\date{}       
\maketitle

\begin{abstract}
\noindent Using the spectral subspaces obtained in \cite{HS}, Brown's results (cf. \cite{Bro}) on the Brown measure of an operator in a type
II$_1$ factor $(\CM,\tau)$ are generalized to finite sets of commuting
operators in $\CM$. It is shown that whenever $T_1,\ldots, T_n\in\CM$ are
mutually commuting operators, there exists one and only one compactly
supported Borel probability measure $\mu_{T_1,\ldots, T_n}$ on $\B(\C^n)$
such that for all $\alpha_1, \ldots, \alpha_n\in\C$,
\[
\tau\big(\log|\alpha_1T_1+\cdots +\alpha_nT_n - \unit|\big)= \int_{\C^n} \log|\alpha_1 z_1
+\cdots+ \alpha_n z_n-1|\,\d\mu_{T_1,\ldots, T_n}(z_1,\ldots,z_n).
\]
Moreover, for every polynomial $q$ in $n$ commuting variables,
$\mu_{q(T_1,\ldots, T_n)}$ is the push-forward measure of $\mu_{T_1,\ldots,T_n}$ via the map $q:\C^n\rightarrow \C$.\\
\indent In addition it is shown that, as in \cite{HS}, for every Borelset $B\subseteq
\C^n$ there is a maximal closed $T_1$-,..., $T_n$-invariant subspace $\CK$
affiliated with $\CM$, such that $\mu_{T_1|_{\CK},\ldots, T_n|_{\CK}}$
is concentrated on $B$. Moreover, $\tau(P_\CK)=\mu_{T_1,\ldots, T_n}(B)$. This generalizes the main result from \cite{HS} to $n$-tuples of commuting operators in $\CM$.
\end{abstract}

\section{Introduction.}

In \cite{Bro} Brown showed that for every operator $T$ in a type II$_1$ factor
$(\CM,\tau)$ there is one and only one compactly supported Borel
probability measure $\mu_T$ on $\C$, the {\it Brown measure} of $T$, such that for all $\lambda\in\C$,
\[
\tau(\log|T-\lambda\unit|)=\int_\C \log|z-\lambda|\,\d\mu_T(z).
\]
In \cite{HS} it was shown that given $T\in\CM$ and $B\in\B(\C)$, there is a
maximal closed $T$-invariant projection $P=P_T(B)\in\CM$, such that the Brown
measure of $PTP$ (considered as an element of $P\CM P$) is concentrated on
$B$. Moreover, $P_T(B)$ is hyperinvariant for $T$,
\begin{equation}
  \tau(P_T(B))=\mu_T(B),
\end{equation}
and the Brown measure of $P^\bot TP^\bot$ (considered as an element of
$P^\bot\CM P^\bot$) is concentrated on $B^c$.

In particular, if $S,T\in\CM$ are commuting operators and $A,B\in\B(\C)$,
then $P_S(A)\wedge P_T(B)$ is $S$- and $T$-invariant. Thus, it is
tempting to define a Brown measure for the pair $(S,T)$, $\mu_{S,T}$, by
\begin{equation}\label{introx}
\mu_{S,T}(A\times B)= \tau(P_S(A)\wedge P_T(B)), \qquad A,B\in\B(\C).
\end{equation}
In order to show that $\mu_{S,T}$ extends (uniquely) to a Borel probability
measure on $\C^2$, we introduce the notion of an {\it idempotent valued
measure} (cf. section~3), and for $T\in\CM$ and $B\in\B(\C)$ we define an
unbounded idempotent affiliated with $\CM$, $e_T(B)$, by
$\CD(e_T(B))=\CK_T(B)+\CK_T(B^c)$ and
\begin{equation}\label{intro1}
e_T(B)\xi = \left\{
    \begin{array}{ll}
      \xi, & \xi \in \CK_T(B)\\
       0, & \xi \in \CK_T(B^c)
    \end{array}\right.
\end{equation}
(cf. section~3). We prove that the map $B\mapsto e_T(B)$ is an idempotent valued
measure and this enables us to show that $\mu_{S,T}$ given by \eqref{introx} has an extension to a
measure on $\C^2$ (cf. section~4). As in \cite{HS} we also prove the
existence of {\it spectral subspaces} for a set of commuting operators
in $\CM$. In the case of two commuting operators $S,T\in\CM$, we show that
for $B\in\B(\C^2)$ there is a maximal $S$- and $T$-invariant projection
$P\in \CP(\CM)$, such that $\mu_{S|_{P(\CH)}, T|_{P(\CH)}}$ (computed
relative to $P\CM P$) is concentrated
on $B$. In the case where $B=B_1\times B_2$, $B_1, B_2\in\B(\C)$, $P$ is
simply the intersection $P_S(B_1)\wedge P_T(B_2)$ (cf. section~5).

Finally, in section~6 we show that $\mu_{S,T}$ has a characterization
similar to the one given by Brown in \cite{Bro}. That is, we show that 
$\mu_{S,T}$ is the unique compactly
supported Borel probability measure on $\B(\C^2)$
such that for all $\alpha, \beta\in\C$,
\[
\tau\big(\log|\alpha S+\beta T - \unit|\big)= \int_{\C^2} \log|\alpha z
+ \beta w -1|\,\d\mu_{S,T}(z,w),
\]
and there is a similar characterization of the Brown measure of an arbitrary finite set of commuting operators. In particular, $\mu_{\alpha S + \beta T}$ is the push-forward measure of
$\mu_{S,T}$ via the map $(z,w)\mapsto \alpha z + \beta w$. In fact we can show that for every polynomial $q$ in two commuting variables, $\mu_{q(S,T)}$ is the push-forward measure of $\mu_{S,T}$ via the map $q:\C^2\rightarrow \C$. All this can (and will) be generalized to arbitrary finite sets of mutually commuting operators. 

Section~2 is devoted to a summary of \cite{Ne} in which we recall the
definition of the measure topology on $\CM$ and how the closure of $\CM$
w.r.t. the measure topology may be identified with the set of closed, densely
defined operators affiliated with $\CM$. Afterwards we focus our attention
on the set of unbounded idempotents affiliated with $\CM$, $\CI(\tilde\CM)$. We prove that these are in
one-to-one correspondence with pairs of projections $P,Q\in \CM$ with
$P\wedge Q =0$ and $P\vee Q=\unit$. More precisely, when $E$ is a closed,
densely defined, unbounded operator affilated with $\CM$ with
$\overline{E\cdot E}=E$, then $P$, the range projection of $E$, and $Q$,
the projection onto the kernel of $E$, satisfy that $P\wedge Q =0$ and $P\vee
Q=\unit$. Moreover, $\CD(E)=P(\CH)+Q(\CH)$, and $E$ is given by
\begin{equation}\label{intro2}
E\xi = \left\{
    \begin{array}{ll}
      \xi, & \xi \in P(\CH)\\
       0, & \xi \in Q(\CH)
    \end{array}\right.
\end{equation}
Vice versa, when $P,Q\in \CM$ with $P\wedge Q =0$ and $P\vee Q=\unit$, then 
\eqref{intro2} defines a closed,
densely defined, unbounded operator affilated with $\CM$ with
$\overline{E\cdot E}=E$. In particular, \eqref{intro1} defines a closed,
densely defined idempotent affilated with $\CM$. In section~2 it is also shown that $\CI(\tilde\CM)$ is stable w.r.t. addition of countably many
idempotents $(E_n)_{n=1}^\infty$ with $E_nE_m=E_mE_n=0$, $n\neq
m$. These are results which are needed in the succeeding sections. 

\vspace{.2cm}

{\it Aknowledgements.} Part of this work was carried out while I visited the UCLA~Department~of {\mbox Mathematics}. I want to thank the department, and especially the Operator~Algebra~Group, for their hospitality. I also thank my advisor, Uffe~Haagerup, with whom I had enlightening discussions about this work. 

\section{Idempotents in $\tilde\CM$.}

We begin this section with a summary of E. Nelson's ``Notes on
non-commutative integration'' \cite{Ne}. We consider a finite von Neumann
algebra $\CM$ represented on a Hilbert space $\CH$ and we fix a faithful,
normal, tracial state $\tau$ on $\CM$. We let
$\CP(\CM)$ denote the set of projections in $\CM$.

Nelson defines the {\it measure topology} on $\CM$ as follows: For $\eps,
\delta>0$, let 
\[
{\rm N}(\eps, \delta)=\{T\in \CM\,|\,\exists P\in \CP(\CM): \|TP\|\leq \eps, \tau(P^\perp)\leq \delta\}.
\]
The measure topology on $\CM$ is then the translation invariant
topology on $\CM$ for which the ${\rm N}(\eps, \delta)$'s form a fundamental
system of neighborhoods of $0$. $\tilde\CM$ is the completion of $\CM$
w.r.t. the measure topology.

Similarly, Nelson defines $\tilde\CH$ to be the completion of $\CH$
w.r.t. the translation invariant topology on $\CH$ for which the sets
\[
{\rm O}(\eps, \delta)=\{\xi\in \CH\,|\,\exists P\in \CP(\CM): \|P\xi\|\leq \eps, \tau(P^\perp)\leq \delta\}
\]
form a fundamental system of neighborhoods of $0$. According to
\cite[Theorem~2]{Ne}, the natural mappings $\CM\rightarrow\tilde\CM$ and
$\CH\rightarrow\tilde\CH$ are both injections.

\begin{thm}\cite[Theorem~1]{Ne} The mappings
\begin{eqnarray*}
\CM\rightarrow \CM&:& T\mapsto T^*,\\
\CM\times \CM \rightarrow \CM&:& (S,T)\mapsto S+T,\\
\CM\times \CM \rightarrow \CM&:& (S,T)\mapsto S\cdot T,\\
\CH\times \CH \rightarrow \CH&:& (\xi,\eta)\mapsto \xi +\eta,\\
\CM\times \CH\rightarrow \CH&:& (T,\xi)\mapsto T\xi
\end{eqnarray*}
all have unique continuous extensions as mappings $\tilde \CM \rightarrow
\tilde \CM$, $\tilde \CM\times \tilde \CM \rightarrow \tilde \CM$, $\tilde
\CM\times \tilde \CM \rightarrow \tilde \CM$, $\tilde \CH\times \tilde \CH
\rightarrow \tilde \CH$ and $\tilde \CM\times \tilde \CH \rightarrow \tilde
\CH$, respectively. In particular, $\tilde \CH$ is a complex vector space
and $\tilde \CM$ is a complex $*$-algebra with a continuous representation
on $\tilde \CH$.
\end{thm}

\vspace{.2cm}

For $x\in\tilde\CM$ define
\begin{equation}\label{corr2}
  \CD(M_x)= \{\xi\in\CH\,|\,x\xi\in\CH\}
\end{equation}
and define $M_x:\CD(M_x)\rightarrow \CH$ by
\begin{equation}
  M_x\xi = x\xi, \qquad \xi\in\CD(M_x).
\end{equation}

Recall that a (not necessarily bounded or everywhere defined) operator $A$
on $\CH$ is said to be {\it affiliated} with $\CM$ if $AU=UA$ for every
unitary $U\in\CM'$. 

\begin{thm}\cite[Theorem~4]{Ne}\label{NeThm4} For every $x\in\tilde\CM$, $M_x$ is a
  closed, densely defined operator affiliated with $\CM$, and
  \begin{equation}
  M_x\cc= M_{x\cc}.
  \end{equation}
  Moreover, for $x,y\in\tilde\CM$,
  \begin{eqnarray}
    M_{x+y}&=&\overline{M_x+M_y},\\
    M_{x\cdot y}&=&\overline{M_x\cdot M_y},
  \end{eqnarray}
  where $\overline{A}$ denotes the closure of a closable operator $A$.
\end{thm}
  
\vspace{.2cm}

A closed, densely defined operator $A$ on $\CH$ has a polar
decomposition $A=V|A|$, and if $A$ is affiliated with $\CM$, then
$V\in\CM$ and all the spectral projections $(E_{|A|}([0,t[))_{t>0}$ belong
to $\CM$. Put
\[
A_n = V\int_0^n t\,\d E_{|A|}(t).
\]
Assuming that $A$ is affiliated with $\CM$, we get that
$(A_n)_{n=1}^\infty$ is a Cauchy sequence w.r.t. the measure
topology. Indeed,
\[
(A_{n+k}-A_n)\cdot E_{|A|}([0,n[)=0,
\]
and $\tau(E_{|A|}([n,\infty[))\rightarrow 0$ as $n\rightarrow
\infty$. Hence, there exists $a\in\tilde\CM$ s.t. $A_n\rightarrow a$ in
measure, and according to \cite[Theorem~3]{Ne}, $A=M_a$. It follows that
every closed, densely defined operator $A$ affiliated with $\CM$ is of the
form $A=M_a$ for some $a\in\tilde\CM$, and this $a$ is uniquely
determined. Indeed, if $M_a=M_b$, then $a$ and $b$ agree on $\CD(M_a)$
which is dense in $\CH$ w.r.t. the norm topology and hence dense in
$\tilde\CH$ w.r.t. the measure topology. Since the representation of
$\tilde\CM$ on $\CH$ is continuous, it follows that $a$ and $b$ agree on
all of $\tilde\CH$. By the same argument, if $S$ and $T$ are
closed, densely defined operators affiliated with $\CM$ which agree on a
dense subset of $\CH$, then $S=T$.

In summary, $\tilde\CM$ is the completion of $\CM$ w.r.t. the measure
topology but it may also be viewed as the set of closed, densely defined
operators affiliated with $\CM$. In particular, if $S$ and $T$ belong to the
latter, then Theorem~\ref{NeThm4} tells us that $S\cc$, $\overline{S+T}$ and
$\overline{S\cdot T}$ are also closed, densely defined and affiliated with
$\CM$.

\begin{notation} In what follows we shall identify $\tilde\CM$ with the set
  of closed, densely defined operators affiliated with $\CM$, but whenever
  necessary, we will specify which one of the two pictures mentioned above we are using. In general,
lower case letters will represent elements of the completion of $\CM$
w.r.t. the measure topology, whereas upper case letters will represent closed,
densely defined operators affiliated with $\CM$. For $x\in\tilde\CM$ we
will denote by $\ker(x)$, $\range(x)$ and $\supp(x)$ the kernel of $M_x$,
the range of $M_x$ and the support projection of $M_x$, respectively.
\end{notation}

\vspace{.2cm}

In the rest of this section we will study the set of idempotent elements in
$\tilde\CM$,
\begin{equation}
\CI(\tilde\CM)=\{e\in\tilde\CM\,|\,e\cdot e=e\}.
\end{equation}
Alternatively, if $\tilde\CM$ is viewed as the set of closed, densely defined operators
affiliated with $\CM$, then $\CI(\tilde\CM)$ is the subset of operators
$E$ fulfilling that $\overline{E\cdot E}=E$. 

The following proposition shows that $\CI(\tilde\CM)$ is in one-to-one
correspondence with the set of pairs of projections $P,Q\in\CM$ such that
$P\wedge Q=0$ and $P\vee Q=\unit$.

\begin{prop}\label{ker-range} Let $E\in \CI(\tilde\CM)$. Then the range of $E$, $\range(E)$,
  and the kernel of $E$, $\ker(E)$ (which is the range of $\unit-E$), are closed subspaces of
  $\CH$, with
\begin{equation}\label{intersect}
  \range(E)\cap \ker(E)=\{0\}
\end{equation}
and
\begin{equation}\label{dense}
  \overline{ \range(E)+\ker(E)}=\CH.
\end{equation}
Moreover, $\CD(E)=\range(E)+ \ker(E)$. Conversely, if $P$, $Q\in\CM$ are projections with $P\wedge Q=0$ and $P\vee
  Q=\unit$, then there is a unique idempotent $E\in\tilde\CM$ with $\CD(E)=P(\CH)+Q(\CH)$, $\range(E)=P(\CH)$ and $\ker(E)=Q(\CH)$, and $E$ is determined by
  \begin{equation}\label{define e}
  E\xi = \left\{
    \begin{array}{ll}
      \xi, & \xi \in P(\CH),\\
       0, & \xi \in Q(\CH).
    \end{array}\right.
  \end{equation}
\end{prop}

\proof Write $E=M_e$ for a (uniquely determined) element $e\in\tilde\CM$
with $e\cdot e=e$ and recall from \eqref{corr2} that
\[
\CD(E)=\{\xi\in\CH\,|\,e\xi\in\CH\}.
\]
Let $E\cdot E$ denote the densely defined operator on $\CH$ obtained by composing $E$ with itself. Then
\begin{eqnarray*}
\CD(E\cdot E)&=& \{\xi\in\CD(E)\,|\, E\xi\in\CD(E)\}\\
&=& \{\xi\in\CD(E)\,|\, eE\xi\in\CH\}\\
&=& \{\xi\in\CD(E)\,|\, e^2\xi\in\CH\}.
\end{eqnarray*}
But $e^2=e$ as everywhere defined operators on $\tilde\CH$. Hence, for $\xi\in\CD(E)$,
\[
e^2\xi = e\xi = E\xi,
\]
and it follows that $\CD(E\cdot E)=\CD(E)$ and that $E\cdot E=E$ (without taking closure). In particular, $\range(E)\subseteq \CD(E)$. Moreover, since $E\cdot E=E$, 
\[
\range(E)= \{\xi\in\CD(E)\,|\,E\xi=\xi\}=\ker(\unit-E).
\]
According to \cite[Exercise~2.8.45]{KR}, the kernel of a closed operator is closed. Hence, $\ker(E)$ and $\range(E)=\ker(\unit-E)$ are closed. Moreover,
\[
\range(E)\cap \ker(\unit-E)= \ker(\unit-E)\cap \ker(E)=\{0\}.
\]
Clearly, $\range(E)+\ker(E)\subseteq \CD(E)$, and since
\[
\xi = E\xi + (\unit-E)\xi, \qquad \xi\in\CD(E),
\]
the converse inclusion also holds. That is,
\[
\CD(E)=\range(E)+\ker(E)
\]
Let $P$ and $Q$ denote the projections onto $\range(E)$ and $\ker(E)$, respectively. Then by the above, $P\wedge Q=0$ and $P\vee Q=\unit$, and $E$ is the operator on $\CD(E)=P(\CH)+Q(\CH)$ determined by \eqref{define e}. In particular, $E$ is uniquely determined by its range and its kernel. 

Conversely, assume that $P$ and $Q$ are projections in $\CM$, such that $P\wedge Q=0$ and $P\vee Q=\unit$, and let $E$ be the operator on $\CD(E)=P(\CH)+Q(\CH)$ determined by \eqref{define e}. Then clearly, $E\cdot E=E$ (without taking closure), $\range(E)=P(\CH)$, and $\ker(E)=Q(\CH)$. Moreover, the graph of $E$ is given by
\begin{eqnarray*}
\CG(E)&=&\{(\xi+\eta,\xi)\,|\,\xi\in P(\CH), \eta\in Q(\CH)\}\\
&=& \{(u,v)\in\CH\times\CH\,|\,u-v\in Q(\CH), v\in P(\CH)\},
\end{eqnarray*}
which is clearly a closed subspace of $\CH\times\CH$. Hence, $E\in\CI(\tilde\CM)$. $\endproof$ 

\vspace{.2cm}

\begin{definition} We define $\tr:\CI(\tilde\CM)\rightarrow [0,1]$ by
  \begin{equation}
    \tr(e)=\tau(\supp(e)), \qquad e\in \CI(\tilde\CM),
  \end{equation}
where $\supp(e)\in \CP(\CM)$ denotes the support projection of $M_e$.
\end{definition}

\vspace{.2cm}

\begin{remark}\label{support-range} For every $x\in\tilde\CM$, $\supp(x)\sim P_{\range(x)}$ so one also has that for $e\in\CI(\tilde \CM)$,
\begin{equation}\label{1-11}
  \tr(e)=\tau(P_{\range(e)}).
\end{equation}
Throughout the paper we will, without further mentioning, make use of this identity. 
\end{remark}

\vspace{.2cm}

\begin{prop}\label{finitesum}
Let $e_1, \ldots, e_n\in\CI(\tilde\CM)$ (seen as everywhere defined operators on $\tilde\CH$) with $e_ie_j=0$ when $i\neq j$. Then
$e_1+\cdots +e_n\in\CI(\tilde\CM)$, and
\begin{itemize}
  \item[(a)]$\ker(e_1+\cdots +e_n)=\bigcap_{i=1}^n \ker(e_i)$,
  \item[(b)]$\supp(e_1+\cdots +e_n)=\bigvee_{i=1}^n \supp(e_i)$,
\item[(c)]$\range(e_1+\cdots +e_n)=\bigvee_{i=1}^n\range(e_i)$,
  \item[(d)]$\tr( e_1+\cdots +e_n)=\sum_{i=1}^n \tr(e_i)$.
\end{itemize} 
\end{prop}

\proof It suffices to consider the case $n=2$. The general case follows by
induction over $n\in\N$. If $e_1, e_2\in\CI(\tilde\CM)$ with
$e_1e_2=e_2e_1=0$, then obviously, $e_1+e_2\in\CI(\tilde\CM)$. 

(a) Clearly, $\ker(e_1)\cap \ker(e_2)\subseteq \ker(e_1+e_2)$. On the other
hand, if $\xi\in \ker(e_1+e_2)$, then $e_i\xi= e_i(e_1+e_2)\xi =0$,
($i=1,2$), whence $\ker(e_1+e_2)\subseteq \ker(e_1)\cap\ker(e_2)$.

(b) Since $\supp(e_1+e_2)$ is the projection onto $[\ker(e_1+e_2)]^\bot$, (b) follows from (a).

(c) For a closed, densely defined operator $S$ affiliated with $\CM$, $\overline{\range(S)}=\ker(S\cc)^\bot$. Using that every idempotent has closed range (cf. Proposition~\ref{ker-range}) and applying (a) to $e_1\cc+\cdots + e_n\cc$, we find that
\begin{eqnarray*}
\range(e_1+\cdots + e_n)&=& \overline{\range(e_1+\cdots + e_n)}\\
&=& [\ker(e_1\cc+\cdots+e_n\cc)]^\bot\\
&=& \bigvee_{i=1}^n \ker(e_i\cc)^\bot\\
&=& \bigvee_{i=1}^n\range(e_i).
\end{eqnarray*}

(d) For $i=1,2$ put $P_i=P_{\range(e_i)}$. Since $e_1e_2=e_2e_1=0$, for $\xi\in P_1(\CH)\cap P_2(\CH)$ we have that
\[
\xi = e_1\xi = e_2e_1\xi = 0.
\]
Then by Kaplansky's formula (cf. \cite[Theorem~6.1.7]{KR}),
\[
P_1\vee P_2 - P_1 \sim P_2 - P_1\wedge P_2 = P_2,
\]
so that
\[
\tau(P_1\vee P_2)= \tau(P_1)+\tau(P_2).
\]
According to (c), $P_1\vee P_2$ is the projection onto $\range(e_1+e_2)$, and then by Remark~\ref{support-range},
\[
\tr(e_1+e_2)=\tr(e_1)+\tr(e_2).\qquad \endproof
\]

\vspace{.2cm}

\begin{prop}\label{infinsum} Let $(e_n)_{n=1}^\infty$ be a sequence of idempotents in
  $\tilde\CM$ with $e_ne_m = e_me_n = 0$ when $n\neq m$. Then
  there is an idempotent in $\tilde\CM$, which we denote by
  $\sum_{n=1}^\infty e_n$, such that
  $\sum_{n=1}^Ne_n\rightarrow \sum_{n=1}^\infty e_n $ in measure as $N\rightarrow\infty$.
  Moreover, $\supp(\sum_{n=1}^Ne_n)\nearrow \supp(\sum_{n=1}^\infty e_n)$,
  whence 
  \begin{equation}\label{Br22}
    \supp\Big(\sum_{n=1}^\infty e_n\Big)=\bigvee_{n=1}^\infty \supp(e_n),
  \end{equation}
  and
  \begin{equation}\label{Br23}
    \tr\Big(\sum_{n=1}^\infty e_n\Big)=\sum_{n=1}^\infty \tr(e_n).
  \end{equation}
  Also,
  \begin{equation}\label{range of sum}
    \range\Big(\sum_{n=1}^\infty e_n\Big)= \bigvee_{n=1}^\infty \range(e_n).
  \end{equation}
\end{prop}

\proof Let
\begin{equation}
  f_n= \sum_{k=1}^ne_k, \qquad n\in\N.
\end{equation}
Then, according to Proposition~\ref{finitesum}, $\supp(f_n)=\bigvee_{k=1}^n
  \supp(e_k)$, and $\tr(f_n)= \sum_{k=1}^n\tr(e_k)$. We prove that
  $(f_n)_{n=1}^\infty$ is a Cauchy sequence in $\tilde\CM$. For $n,
  k\in\N$, let
  \begin{equation}
    P_{n,k}= \supp(f_{n+k}-f_n)= \bigvee_{l=1}^k\supp(e_{n+l}).
  \end{equation}
  Then
  \begin{equation}
    (f_{n+k}-f_n)P_{n,k}^\bot =0,
  \end{equation}
  so for every $\eps>0$,
  \begin{equation}\label{Br10}
    f_{n+k}-f_n\in {\rm N}(\eps, \tau(P_{n,k}))={\rm N}\Big(\eps, \sum_{l=1}^k \tr(e_{n+l})\Big).
  \end{equation}
  Now, $\sum_{k=1}^n \tr(e_{k}) =\tr(f_n)\leq 1$, so
  for arbitrary $\delta>0$ there is an $n_0\in\N$ such that
  \begin{equation}\label{Br11}
    \sum_{k=n_0}^\infty \tr(e_k)\leq \delta.
  \end{equation}
  It follows from \eqref{Br10} and \eqref{Br11} that when $n\geq n_0$ and
  $k\geq 1$, then $f_{n+k}-f_n\in {\rm N}(\eps, \delta)$. Thus,
  $(f_n)_{n=1}^\infty$ is a Cauchy sequence in $\tilde\CM$. Put
  \begin{equation}
    e= \lim_{n\rightarrow\infty}f_n\in \tilde\CM.
  \end{equation}
For all $k,n\in\N$, $f_{n+k}f_n=f_nf_{n+k}=f_n$, and hence
  \begin{equation}\label{Br12}
    ef_n=f_n e=f_n,
  \end{equation}
  so that $e\cdot e=e$.
  
  Let $P_n=\supp(f_n)=\bigvee_{k=1}^n \supp(e_k)$. Then $P_n\nearrow
  P:=\bigvee_{k=1}^\infty \supp(e_k)$ as $n\rightarrow \infty$, and
  \begin{eqnarray*}
    \sum_{k=1}^\infty\tr(e_k) &=&
    \lim_{n\rightarrow\infty}\sum_{k=1}^n\tr(e_k)\\
     &=& \lim_{n\rightarrow\infty}\tr(f_n)\\
     &=&\lim_{n\rightarrow\infty}\tau(P_n)\\
     &=&\tr(P).
  \end{eqnarray*}
  It follows from \eqref{Br12} that for every $n\in\N$, $P_n\leq
  \supp(e)$. Hence $P\leq \supp(e)$. On the other hand, for every $n\in\N$,
  $f_n(\unit-P_n)=0$, so
  \begin{equation}
    e(\unit-P)=\lim_{n\rightarrow\infty}[f_n(\unit-P_n)]=0
  \end{equation}
  (the limit refers to the measure topology). Thus, $\supp(e)\leq P$, and
  we have shown that $\supp(e)=P=\bigvee_{k=1}^\infty \supp(e_k)$ and that
  \eqref{Br23} holds.

In order to prove \eqref{range of sum}, let $\xi\in \range(e_m)$. Then
$e_m\xi = \xi$ and
\[
\Big(\sum_{n=1}^\infty e_n\Big)\xi = e_m\xi = \xi.
\]
Thus,
\begin{equation}
\range\Big(\sum_{n=1}^\infty e_n\Big)\supseteq \bigvee_{n=1}^\infty
\range(e_n).
\end{equation}
On the other hand, we know that $\range\Big(\sum_{n=1}^\infty
e_n\Big)\subseteq \CH$, and since $\range\Big(\sum_{n=1}^N e_n\Big)
\subseteq  \bigvee_{n=1}^\infty \range(e_n)$ for all $N\in\N$, we also have
that
\[
\range\Big(\sum_{n=1}^\infty e_n\Big)\subseteq {\Big(\bigvee_{n=1}^\infty
\range(e_n)\Big)}\tilde{},
\]
where $\tilde{}$ denotes closure w.r.t. the measure topology. Intersecting by
$\CH$ on both sides of the inclusion, we get that 
\begin{equation}
\range\Big(\sum_{n=1}^\infty e_n\Big)\subseteq \bigvee_{n=1}^\infty
\range(e_n).
\end{equation}
This proves \eqref{range of sum}. $\endproof$

\vspace{.2cm}

We shall make use of the following theorem from \cite{A}. For a published
proof of it, we refer the reader to \cite{Aa}. 

\begin{thm}\cite{A} \label{Ainsworth} Let $E$ and $F$ be (not necessarily closed)
  subspaces of $\CH$ which are affiliated with $\CM$.\footnote{A subspace
  $E$ of $\CH$ is said to be {\it affiliated with $\CM$} if for all
  $T\in\CM'$, $T(E)\subseteq E$. Note that if $E$ is affiliated with $\CM$,
  then the projection onto $\overline{E}$ belongs to $\CM$, and if $E$ is
  {\it closed}, then this is a necessary and sufficient condition for $E$ to be
  affiliated with $\CM$.} Then $E\cap F$ is affiliated with $\CA$, and
\begin{equation}
  \overline{E\cap F}=\overline{E}\cap\overline{F}.
\end{equation}
\end{thm}

  \vspace{.2cm}
  
  \begin{lemma} Consider idempotents $e$, $f\in\tilde\CM$. Let $P=P_{\range(e)}$,
  $Q=P_{\range(\unit-e)}$, $R=P_{\range(f)}$ and $S=P_{\range(\unit-f)}$. Then $ef=fe$ if and only if
  \begin{equation}\label{Br9}
    (P\wedge R)\vee(P\wedge S)\vee (Q\wedge R)\vee (Q\wedge S)=\unit.
  \end{equation}
\end{lemma}

\proof Clearly,
\begin{equation}
  \unit = ef + e(\unit-f) + (\unit-e)f + (\unit-e)(\unit-f).
\end{equation}
Suppose that $ef=fe$, and let $g_1= ef$, $g_2=e(\unit-f)$, $g_3= (\unit-e)f$ and $g_4=
(\unit-e)(\unit-f).$ Then $g_1,\ldots, g_4$ are idempotents with support
projections $P_1, P_2, P_3$ and $P_4$, respectively, such that
\begin{equation}
  \big(\bigvee_{i=1}^4P_i\big)(\CH)=\CH.
\end{equation}
Moreover, $P_1\leq P\wedge R$, $P_2\leq P\wedge S$, $P_3\leq Q\wedge R$ and
$P_4\leq Q\wedge S$. This shows that \eqref{Br9} holds.

  On the other hand, assume that  \eqref{Br9} holds. According to
  \eqref{Br9} and Theorem~\ref{Ainsworth},
  \[
  \CH_0 := \CD(ef)\cap\CD(fe)\cap  [(P\wedge R)(\CH)+(P\wedge S)(\CH)+ (Q\wedge R)(\CH)+
  (Q\wedge S)(\CH)]
  \]
is dense in $\CH$ so it suffices to prove that $ef$ and $fe$ agree on
$\CH_0$. To see this, let $\xi\in\CD(ef)\cap\CD(fe)\cap(P\wedge R)(\CH)$. Then
  \begin{equation}
    ef\xi = e\xi = \xi = f\xi = fe\xi ,
  \end{equation}
  and similarly, when $\xi \in \CD(ef)\cap\CD(fe)\cap (P\wedge S)(\CH)$, $\xi \in
  \CD(ef)\cap\CD(fe)\cap (Q\wedge R)(\CH)$ or  $\xi \in
  \CD(ef)\cap\CD(fe)\cap (Q\wedge S)(\CH)$. Thus, $ef$ agrees with $fe$ on
  $\CH_0$. $\endproof$ 

\section{An idempotent valued measure associated with $T\in\CM$.}

As in the previous section, consider a finite von Neumann algebra $\CM$
with a faithful, normal, tracial state $\tau$. Inspired by the notion of a spectral measure we make the following defintion:

\begin{definition}\label{skewprojmeas} Let $(X,\CF)$ denote a measurable space. An {\it
    idempotent valued measure} on $(X,\CF)$ (with values in $\tilde\CM$) is a map $e$ from $\CF$
    into $\CI(\tilde\CM)$ such that
    \begin{itemize}
    \item[(i)]$e(X)=\unit$,
    \item[(ii)]  $e(F_1)e(F_2)=e(F_2)e(F_1)=0 $ when $F_1, F_2\in \CF$ with
    $F_1\cap F_2=\emptyset$,
    \item[(iii)] when $(F_n)_{n=1}^\infty$ is a sequence of mutually
    disjoint sets from $\CF$, then $\sum_{n=1}^N e(F_n)$ converges in
    measure as $N\rightarrow\infty$ to $e\Big(\bigcup_{n=1}^\infty
    F_n\Big)$, i.e.
    \begin{equation*}
      e\Big(\bigcup_{n=1}^\infty F_n\Big)=\sum_{n=1}^\infty e(F_n).
    \end{equation*}
    \end{itemize}
\end{definition}

Note that because of (ii) and Proposition~\ref{infinsum}, the limit in
(iii) actually exists. 

\vspace{.2cm}

From now on we will assume that $\CM$ is in fact a type II$_1$ factor. Recall from \cite{HS} that for $T\in\CM$ and $B\subseteq \C$ a Borel set there is a
maximal $T$-invariant projection $P=P_T(B)\in\CM$, such that the Brown
measure of $PTP$ (considered as an element of $P\CM P$) is concentrated on
$B$. Moreover, $P_T(B)$ is hyperinvariant for $T$,
\begin{equation}\label{Br1}
  \tau(P_T(B))=\mu_T(B),
\end{equation}
and the Brown measure of $P^\bot TP^\bot$ (considered as an element of
$P^\bot\CM P^\bot$) is concentrated on $B^c$. We let $\CK_T(B)$ denote the
range of $P_T(B)$. Then the aim of this section is to
prove:

\begin{thm}\label{skewprojmeas-T} Let $T\in\CM$, and for $B\in \B(\C)$, let
  $e_T(B)$ with $\CD(e_T(B))=\CK_T(B)+\CK_T(B^c)$ be given by
  \begin{equation}\label{Br8}
    e_T(B)\xi = \left\{
    \begin{array}{ll}
      \xi, & \xi \in \CK_T(B)\\
       0, & \xi \in \CK_T(B^c)
    \end{array}\right.
\end{equation}
Then $e_T(B)\in\CI(\tilde \CM)$, and $B\mapsto e_T(B)$ is an idempotent valued measure.
\end{thm}

\vspace{.2cm}

The proof of this theorem uses various results which we state and prove
below. The first one of these is a lemma which we proved in
\cite{HS}, but for the sake of completeness we give the proof here as well.

\begin{lemma}\label{intersect-inv} Let $T\in\CM$, and let $P\in\CM$ be a
  non-zero, $T$-invariant projection. Then for every $B\in\B(\C)$,
  \begin{equation}
    \CK_{T|_{P(\CH)}}(B)=\CK_T(B)\cap P(\CH),
  \end{equation}
  where $T|_{P(\CH)}$ is considered as an element of the type II$_1$ factor $P\CM P$.
\end{lemma}

\proof Let $Q\in P\CM P$ denote the projection onto $ \CK_{T|_{P(\CH)}}(B)$, and
let $R = P_T(B)\wedge P$. We will prove that $Q\leq R$ and $R\leq Q$.

Clearly, $Q\leq P$. In order to see that $Q\leq P_T(B)$, recall that
$P_T(B)$ is maximal w.r.t. the properties
\begin{itemize}
  \item[(i)] $P_T(B)TP_T(B)=TP_T(B)$,
  \item[(ii)] $\mu_{P_T(B)TP_T(B)}$ (computed relative to $P_T(B)\CM P_T(B)$) is concentrated on $B$.
\end{itemize}

Since
\begin{equation}
  QTQ = QTPQ = TPQ = TQ,
\end{equation}
and $\mu_{QTQ}= \mu_{QTPQ}$ (computed relative to $Q\CM Q$) is concentrated
on $B$, we get that $Q\leq P_T(B)$, and hence $Q\leq R$.

Similarly, to prove that $R\leq Q$, we must show that
\begin{itemize}
  \item[(i')] $RTPR = TPR$, i.e. $RTR = TR$,
  \item[(ii')] $\mu_{RTPR}= \mu_{RTR}$  (computed relative to $R\CM R$) is concentrated on $B$.
\end{itemize}

Note that if $P_T(B)=0$, then $R\leq Q$, so we may assume that $P_T(B)\neq 0$.
(i') holds, because $R(\CH)= P(\CH)\cap P_T(B)(\CH)$ is $T$-invariant when
$P(\CH)$ and $P_T(B)(\CH)$ are $T$-invariant. In order to prove (ii'), at
first note that $R(\CH)$ is $TP_T(B)$-invariant. Hence
\begin{equation}
  \mu_{TP_T(B)}= \tau_1(R)\cdot \mu_{RTR} + \tau_1(R^\bot)\cdot
  \mu_{R^\bot TR^\bot},
\end{equation}
where $\tau_1=\frac{1}{\tau(P_T(B))}\cdot\tau|_{P_T(B)\CM P_T(B)}$. It follows that
\begin{equation}
  \tau_1(R)\cdot\mu_{RTR}(B^c)\leq \mu_{TP_T(B)}(B^c)=0,
\end{equation}
and thus, if $R\neq 0$, then $\mu_{RTR}(B^c)=0$, and (ii') holds. If $R=0$,
then $R\leq Q$ is trivially fulfilled. $\endproof$

\vspace{.2cm}

\begin{prop}\label{intersect-complement} For every Borel set $B\subseteq
  \C$,
\begin{equation}\label{Br18}
  \CK_T(B)=\CK_{T\cc}((B^c)\cc)^\bot,
\end{equation}
where $A\cc := \{\overline z\,|\,z\in A\}$ for $A\subseteq \C$. Moreover, for all Borel sets $A, B\subseteq \C$,
\begin{equation}\label{Br2}
  \CK_T(A)\cap \CK_T(B)=\CK_T( A\cap B),
\end{equation}
and
\begin{equation}
 \CK_T( A\cup B)=\overline{\CK_T(A)+\CK_T( B)}.
\end{equation}
\end{prop}

\proof Let $B\in\B(\C)$
and let $P=P_T(B)$. Then $P^\bot$ is $T\cc$-invariant, and 
\begin{equation}
  \mu_{P^\bot T\cc P^\bot}(B\cc)= \mu_{(P^\bot T\cc
  P^\bot)\cc}(B)=\mu_{P^\bot T P^\bot}(B)=0
\end{equation}
(recall that $\mu_{P^\bot T P^\bot}$ is concentrated on $B^c$). Thus, $
\mu_{P^\bot T\cc P^\bot}$ is concentrated on $\C\setminus B\cc$, and
maximality of $P_{T\cc}( \C\setminus B\cc)$ implies that
\begin{equation}\label{Br17}
  P_T(B)^\bot = P^\bot \leq P_{T\cc}( \C\setminus B\cc).
\end{equation}
Since
\begin{eqnarray*}
  \tau( P_{T\cc}( \C\setminus B\cc))&=& \mu_{T\cc}(\C\setminus B\cc)\\
  &=& 1-\mu_{T\cc}(B\cc)\\
  &=& 1-\mu_T(B)\\
  &=& \tau( P_T(B)^\bot),
\end{eqnarray*}
we get from \eqref{Br17} that $P_T(B)^\bot = P_{T\cc}( \C\setminus B\cc)$.

Next, let $A,B\in\B(\C)$. By maximality of $P_T(A)$ and $P_T(B)$, $P_T(
A\cap B)\leq P_T(A)\wedge P_T( B)$, so $\supseteq$ holds in \eqref{Br2}. We let
$\CK :=  \CK_T(A)\cap \CK_T(B)$, and we let $Q$ denote the projection
onto $\CK$. Then, according to Lemma~\ref{intersect-inv},
\begin{eqnarray*}
  \CK &=& \CK_{T|_{\CK_T(A)}}(B)\\
    &=& \CK_{T|_{\CK_T(B)}}(A),
\end{eqnarray*}
proving that $\mu_{QTQ}$ is concentrated on $A$ and on $B$ and therefore on
$A\cap B$. Consequently, $Q\leq P_T( A\cap B)$, so $\subseteq$ also
holds in \eqref{Br2}.

Finally, we infer from \eqref{Br18} and \eqref{Br2} that
\begin{eqnarray*}
  \CK_T( A\cup B)&=& \CK_T( (A^c\cap B^c)^c)\\
  &=& \CK_{T\cc}( (A^c\cap B^c)\cc)^\bot\\
  &=& \CK_{T\cc}( (A^c)\cc\cap (B^c)\cc)^\bot\\
  &=& [\CK_{T\cc}((A^c)\cc)\cap\CK_{T\cc}( (B^c)\cc)]^\bot\\
  &=& \overline{\CK_{T\cc}((A^c)\cc)^\bot + \CK_{T\cc}( (B^c)\cc)^\bot}\\
  &=& \overline{\CK_T( A)+\CK_T( B)}. \quad \endproof
\end{eqnarray*}

\vspace{.2cm}

It follows from Proposition~\ref{intersect-complement} and
Proposition~\ref{ker-range} that for
$B\in\B(\C)$, $e_T(B)$ given by \eqref{Br8} belongs to $\CI(\tilde \CM)$,
as stated in Theorem~\ref{skewprojmeas-T}.

\begin{lemma}\label{tendtozero}Let $(x_n)_{n=1}^\infty$ be a sequence in
  $\tilde\CM$, and suppose $\tau(\supp(x_n))\rightarrow 0$ as $n\rightarrow
  \infty$. Then $x_n\rightarrow 0$ in the measure topology.
\end{lemma}

\proof This is standard. $\endproof$

\vspace{.2cm}

If $S\in\CM$ commutes with $T\in\CM$, then for every $B\in \B(\C)$, $\CK_T(B)$
and $\CK_T(B^c)$ are $S$-invariant, and therefore $S$ commutes with
$e_T(B)$ as well. We prove that, as a consequence of this,
$[e_S(A),e_T(B)]=0$ for every $A\in \B(\C)$

\begin{lemma}\label{[e,T]} Let $T\in\CM$, and let $e\in\CI(\tilde\CM)$ with
  $[e,T]=0$. Then for every $B\in\B(\C)$, $[e,e_T(B)]=0$. In particular, if
  $S\in\CM$ commutes with $T$, then $[e_S(\cdot),e_T(\cdot)]=0 $.
\end{lemma}

\proof Let $P=P_{\range(e)}$, $Q=P_{\range(\unit-e)}$, $R=P_{\range(e_T(B))}$ and
$S=P_{\range(\unit-e_T(B))}$. We prove that \eqref{Br9} holds. Since $eT=Te$,
$P(\CH)$ and $Q(\CH)$ are $T$-invariant. Then by Lemma~\ref{intersect-inv},
\begin{eqnarray*}
  \CK_{T|_{P(\CH)}}(B)&=& \CK_T(B)\cap P(\CH)\\
  &=& R(\CH)\cap P(\CH),
\end{eqnarray*}
and
\begin{eqnarray*}
   \CK_{T|_{P(\CH)}}(B^c)&=& \CK_T(B^c)\cap P(\CH)\\
  &=& S(\CH)\cap P(\CH).
\end{eqnarray*}
Hence $(R\wedge P)\vee (S\wedge P)=P$, and similarly,  $(R\wedge Q)\vee
(S\wedge Q)=Q$. It follows that
\[
\unit = P\vee Q=(R\wedge P)\vee (S\wedge P)\vee  (R\wedge Q)\vee
(S\wedge Q),
\]
as desired. $\endproof$

\vspace{.2cm}

{\it Proof of Theorem~\ref{skewprojmeas-T}.} $e_T(\emptyset)=0$, because
$P_T(\emptyset)=0$. If $B_1, B_2\in \B(\C)$ with
$B_1\cap B_2=\emptyset$, then $\CK_T(B_1)\cap \CK_T(B_2)=\{0\}$,
i.e. $\range(e_T(B_1))\cap \range(e_T(B_2))=\{0\}$. According to
Lemma~\ref{[e,T]}, $[e_T(B_1),e_T(B_2)]=0$ so that
$e_T(B_1)e_T(B_2)\in\CI(\tilde\CM)$. Moreover,
\[
\range(e_T(B_1)e_T(B_2))\subseteq \range(e_T(B_1))\cap
\range(e_T(B_2))=\{0\},
\]
and we conclude that $e_T(B_1)e_T(B_2)=e_T(B_2)e_T(B_1)=0$.

Now, let $(B_n)_{n=1}^\infty$ be a sequence of mutually disjoint Borel
sets. Then for each $N\in\N$ we get from
Proposition~\ref{intersect-complement} and Lemma~\ref{finitesum} that
\begin{eqnarray*}
  \range\Big(e_T\Big(\bigcup_{n=1}^NB_n\Big)\Big)&=&\CK_T\Big(\bigcup_{n=1}^NB_n\Big)\\
  &=& \overline{\CK_T(B_1)+\cdots +\CK_T(B_N)}\\
  &=& \range(e_T(B_1)+\cdots + e_T(B_N))
\end{eqnarray*}
and
\begin{eqnarray*}
  \ker\Big(e_T\Big(\bigcup_{n=1}^NB_n\Big)\Big)&=&
  \CK_T\Big(\Big(\bigcup_{n=1}^NB_n\Big)^c\Big)\\
  &=& \CK_T\Big(\bigcap_{n=1}^NB_n^c\Big)\\
  &=& \bigcap_{n=1}^N\CK_T(B_n^c)\\
  &=& \bigcap_{n=1}^N \ker(e_T(B_n))\\
  &=& \ker(e_T(B_1)+\cdots + e_T(B_N)).
\end{eqnarray*}
Since an element $e$ in $\CI(\tilde\CM)$ is uniquely determined by its
kernel and its range, it follows that $e_T$ is additive, i.e.
\begin{equation}
  e_T\Big(\bigcup_{n=1}^NB_n\Big)= e_T(B_1)+\cdots + e_T(B_N), \qquad (N\in\N).
\end{equation}

Additivity of $e_T$ implies that
\begin{eqnarray}\label{Br15}
 e_T\Big(\bigcup_{n=1}^\infty B_n\Big)-\sum_{n=1}^\infty
 e_T(B_n)&=&\lim_{N\rightarrow \infty}\Big( e_T\Big(\bigcup_{n=1}^\infty
 B_n\Big)-\sum_{n=1}^N e_T(B_n))\Big)\nonumber \\
 &=& \lim_{N\rightarrow \infty} e_T\Big(\bigcup_{n=N+1}^\infty B_n\Big)
\end{eqnarray}
(the limits refer to the measure topology), where
\begin{eqnarray*}
\tau\Big(\supp\Big(e_T\Big(\bigcup_{n=N+1}^\infty B_n\Big)\Big)\Big)& =&
\tau\Big(P_T\Big(\bigcup_{n=N+1}^\infty B_n\Big)\Big)\\
& = & \mu_T\Big(\bigcup_{n=N+1}^\infty B_n\Big)\\
&\rightarrow & 0, \quad {\rm as} \; N\rightarrow\infty.
\end{eqnarray*}
Combining this with \eqref{Br15} and Lemma~\ref{tendtozero}, we find that
$e_T$ is $\sigma$-additive as well. $\endproof$
  
\vspace{.2cm}

Note that in the case where $T$ is a {\it normal} operator, $B\mapsto
P_T(B)$ is just the spectral measure of $T$, and $e_T(B)=P_T(B)$.

\vspace{.2cm}

\section{The Brown measure of a set of commuting operators in $\CM$.}

As in the previous section, let $\CM$ be a type II$_1$ factor. The purpose of
this section is to prove:

\begin{thm}\label{commuting} Let $n\in\N$, and let $T_1,\ldots, T_n\in \CM$ be commuting operators. Then there is a
  probability measure $\mu_{T_1,\ldots, T_n}$ on $(\C^n, \B(\C^n))$, which is uniquely
  determined by
  \begin{equation}\label{Br4}
    \mu_{T_1,\ldots, T_n}(B_1\times\cdots\times B_n)=
    \tau\Big(\bigwedge_{i=1}^n P_{T_i}(B_i)\Big), \qquad B_1,\ldots, B_n\in
    \B(\C),
  \end{equation}
  where $P_{T_i}(B_i)\in \CM$ is the projection onto $\CK_{T_i}(B_i)$
  (cf. Section~2).
\end{thm}
\vspace{.2cm}

The idea of proof is as follows:

As mentioned in the previous section (cf. Lemma~\ref{[e,T]}), if $S\in\CM$ commutes with $T\in\CM$, then
$[e_S(A),e_T(B)]=0$ for all $A, B\in \B(\C)$. We may therefore define a map $e_{T_1,\ldots, T_n}$ from
$\B(\C)^n$ into $\CI(\tilde\CM)$ by
\begin{equation}\label{Br5}
  e_{T_1,\ldots, T_n}(B_1, \ldots, B_n)= e_{T_1}(B_1) e_{T_2}(B_2)\cdots e_{T_n}(B_n) , \qquad B_1,\ldots, B_n\in\B(\C).
\end{equation}

We will then define $\nu$ on $\B(\C)^n$ by
\begin{equation}
  \begin{split}
   \nu(B_1,\ldots, B_n)&= \tau(\supp[ e_{T_1,\ldots, T_n}(B_1, \ldots,
   B_n)]) \\
   &= \tau(P_{\range(e_{T_1,\ldots, T_n}(B_1, \ldots,
   B_n))})  \\ 
   &=\tau\Big(\bigwedge_{i=1}^n P_{T_i}(B_i)\Big), \qquad B_1,\ldots,
   B_n\in\B(\C), \label{Br6}
   \end{split}
\end{equation}
and we will prove that $\nu$ extends (uniquely) to a probability measure,
$\mu_{T_1,\ldots, T_n}$, on $(\C^n, \B(\C^n))$. 

\vspace{.2cm}

\begin{thm}\label{fundamental} Consider uncountable, complete, separable metric spaces
  $(X_1, d_1), \ldots, (X_n, d_n)$. Suppose $\nu: \B(X_1)\times \cdots\times \B(X_n)\rightarrow [0,\infty[$ is a map satisfying
  \begin{itemize}
    \item[(1)] for all $B_2\in \B(X_2),B_3\in\B(X_3), \ldots, B_n\in \B(X_n)$, $B\mapsto
    \nu(B,B_2, \ldots, B_n)$ is a measure on
    $(X_1, \B(X_1))$,
    \item[(2)] for all $B_1\in \B(X_1),B_3\in\B(X_3), \ldots, B_n\in
    \B(X_n)$, $B\mapsto\nu(B_1, B, B_3,\ldots, B_n) $ is a measure on
    $(X_2, \B(X_2))$,
    
    $\vdots$

    \item[(n)] for all $B_1\in \B(X_1),B_2\in\B(X_2), \ldots, B_{n-1}\in
    \B(X_{n-1})$, $B\mapsto\nu(B_1, B_2, \ldots, B_{n-1}, B) $ is a measure on
    $(X_n, \B(X_n))$.
  \end{itemize}
Then there is a unique measure $\mu$ on $\otimes_{i=1}^n \B(X_i)$, such that for
all $B_1\in \B(X_1),B_2\in\B(X_2), \ldots, B_n\in \B(X_n)$,
\begin{equation}\label{Br7}
  \mu(B_1\times B_2\times \cdots \times B_n)= \nu(B_1, B_2, \ldots, B_n).
\end{equation}
\end{thm}

\proof According to \cite[Remark~1, p.~358]{Ku}, $(X_i, \B(X_i))$ is
Borel equivalent to $([0,1], \B([0,1]))$, i.e. there is a bijective
bimeasurable map $\phi_i: (X_i,\B(X_i))\rightarrow ([0,1], \B([0,1]))$. Therefore we may as well assume
that $X_i=\R$, $(i=1, \ldots, n)$. We may also assume that $\nu(\R,\R,
\ldots, \R)=1$.  We define $F:\R^n \rightarrow [0,1]$ by
\begin{equation}
  F(x_1,\ldots, x_n)=\nu(]-\infty, x_1],\ldots, ]-\infty,x_n]), \qquad
  x_1, \ldots, x_n\in\R.
\end{equation}
Because of (1)--(n), $F$ is increasing in each variable separately and satisfies
\begin{itemize}
  \item[(a)] if $x_i^{(k)}\searrow x_i$, $i=1, \ldots, n$, then
  $F(x_1^{(k)}, \ldots, x_n^{(k)})\searrow F(x_1,\ldots, x_n)$,
  \item[(b)] if $x_i\searrow -\infty$ for some $i\in\{1, \ldots, n\}$, then $F(x_1,\ldots, x_n)
  \searrow 0$,
  \item[(c)] if $x_i\nearrow\infty$ for all $i\in\{1, \ldots, n\}$, then $F(x_1,\ldots, x_n)\nearrow 1$.
\end{itemize}
Then, according to \cite[Corollary~2.27]{Bre}, there is a (unique)
probability measure $\mu$ on $(\R^n, \B(\R^n))$ such that for all
$x_1,\ldots, x_n\in\R$,
\begin{equation}
  \mu(]-\infty, x_1]\times\cdots \times ]-\infty, x_n])=F(x_1,\ldots, x_n).
\end{equation}

Let $x_2, \ldots, x_n\in\R$ be fixed but arbitrary. Then the (finite)
measures
\[
B\mapsto
\mu(B\times ]-\infty, x_2]\times\cdots\times ]-\infty,x_n])
\]
and
\[
B\mapsto \nu(B,]-\infty, x_2], \ldots, ]-\infty,x_n])
\]
have the same
distribution functions. Hence they must be identical. That is, for all
$B\in \B(\R)$,
\begin{equation}\label{Br16}
  \mu(B\times ]-\infty, x_2]\times\cdots\times ]-\infty,x_n] )=\nu(B,]-\infty, x_2], \ldots, ]-\infty,x_n] ).
\end{equation}
Now, let $B_1\in\B(\R)$ and $x_3, \ldots, x_n\in\R$ be fixed but arbitrary. Then \eqref{Br16} shows that
the (finite) measures
\[
B\mapsto \mu(B_1\times B\times]-\infty,
x_3]\times\cdots\times ]-\infty,x_n] )
\]
and
\[
B\mapsto \nu(B_1, B, ]-\infty,
x_3], \ldots, ]-\infty,x_n] )
\]
have 
the same distribution functions, so they must be identical as well. That
is, for all
$B\in \B(\R)$,
\begin{equation}
  \mu(B_1\times B \times ]-\infty, x_3]\times\cdots\times ]-\infty,x_n]) =\nu(B_1,B, ]-\infty, x_3], \ldots, ]-\infty,x_n] ).
\end{equation}
Continuing like this we find that \eqref{Br7} holds. $\endproof$

\vspace{.2cm}

It follows from Theorem~\ref{fundamental} that in order to show that
$\mu_{T_1, \ldots, T_n}$ exists, we must prove that (1)--(n) of
Theorem~\ref{fundamental} hold in the
case where $X_1=\cdots =X_n=\C$, and where $\nu$ is given by \eqref{Br6}.

\vspace{.2cm}

From now on we will, in order to simplify notation a little, consider the case $n=2$,
and we will assume that $S,T\in\CM$ are commuting operators. It should be
clear that the proof given below may be generalized to the case
of arbitrary $n\in\N$. 

\begin{lemma}\label{conditional} For fixed $A\in\B(\C)$, $\nu_A:
  \B(\C)\rightarrow [0,1]$ given by
  \begin{equation}
    \nu_A(B)=\nu(A,B)=\tau(P_S(A)\wedge P_T(B)) , \qquad B\in\B(\C)
  \end{equation}
  (cf. \eqref{Br6}) is a measure on $(\C, \B(\C))$.
\end{lemma}

\proof According to Theorem~\ref{skewprojmeas-T} and Definition~\ref{skewprojmeas},
$e_T(\emptyset)=0$, so
\[
\nu_A(\emptyset)=\tau(\supp[e_S(A)e_T(\emptyset)])=0.
\]
Let $(B_n)_{n=1}^\infty$ be a sequence of mutually disjoint sets from
$\B(\C)$. Then $e_T\Big( \bigcup_{n=1}^\infty B_n\Big) = \sum_{n=1}^\infty
e_T( B_n)$, so
\begin{equation}
  e_S(A)e_T\Big(\bigcup_{n=1}^\infty B_n\Big)= \sum_{n=1}^\infty e_S(A)e_T(
  B_n)
\end{equation}
with $e_S(A)e_T(B_n)e_S(A)e_T(B_m)=e_S(A)e_T(B_n)e_T(B_m)=0$ when $n\neq
m$. Hence, by Proposition~\ref{infinsum},
\begin{equation}
\tr\Big( e_S(A)e_T\Big(\bigcup_{n=1}^\infty B_n\Big)\Big)=\sum_{n=1}^\infty \tr(e_S(A)e_T(B_n)).
\end{equation} 
This shows that $\nu_A$ is a measure. $\endproof$

\vspace{.2cm}

It now follows from Lemma~\ref{conditional} and Theorem~\ref{fundamental} that there
is one and only one (probability) measure $\mu_{S,T}$ on $\B(\C^2)$ such
that for all $A,$ $B\in\B(\C)$,
\begin{equation}
  \mu_{S,T}(A\times B) =  \tau(\supp[e_{S,T}(A,B)])=\tau(P_S(A)\wedge
  P_T(B)),
\end{equation}
and this proves Theorem~\ref{commuting} in the case $n=2$.

\section{Spectral subspaces for commuting operators $S,T\in\CM$.}

\begin{thm}\label{spectral subspaces} Let $S,$ $T\in\CM$ be commuting
  operators, and let $B\subseteq \C^2$ be any Borel set. Then there is a maximal, closed, $S$- and
  $T$-invariant subspace $\CK = \CK_{S,T}(B)$ affiliated with $\CM$, such
  that the Brown measure $\mu_{S|_\CK, T|_\CK}$ is concentrated on
  $B$. Let $P_{S,T}(B)\in\CM$ denote the projection onto
  $\CK_{S,T}(B)$. Then more precisely,
  \begin{itemize}
    \item[(i)] if $B=B_1\times B_2$ with $B_1, B_2\in\B(\C)$, then
  \begin{equation}\label{Br21}
    P_{S,T}(B)=P_S(B_1)\wedge P_T(B_2),
  \end{equation}
    \item[(ii)] if $B$ is a disjoint union of the sets $ (B_1^{(k)}\times
    B_2^{(k)})_{k=1}^\infty$, where $B_i^{(k)}\in\B(\C)$, $k\in\N$, $i=1,2$, then
    \begin{equation}\label{Br19}
      P_{S,T}(B)=\bigvee_{k=1}^\infty [P_S(B_1^{(k)})\wedge
      P_T(B_2^{(k)})],
    \end{equation}
    \item[(iii)] and for general $B\in\B(\C^2)$,
      \begin{equation}\label{Br20}
      P_{S,T}(B)=\bigwedge_{B\subseteq U,\;U\subseteq \C^2 \;{\rm open}}P_{S,T}(U).
      \end{equation}
  \end{itemize}
  Moreover,
  \begin{equation}\label{measure-trace}
    \mu_{S,T}(B)=\tau(P_{S,T}(B)), \quad B\in \B(\C^2).
  \end{equation}
\end{thm}

\vspace{.2cm}

\begin{remark}\label{open sets} Every non-empty, open subset of $\C^2\cong \R^4$
  is a disjoint union of countably many {\it standard intervals}, i.e. sets of
  the form $\prod_{i=1}^4]a_i,b_i]$, where $-\infty<a_i<b_i<\infty$,
  $1\leq i\leq 4$. Hence, it follows from Theorem~\ref{spectral
  subspaces} that the map $\B(\C^2)\rightarrow \CP(\CM): B\mapsto
  P_{S,T}(B)$ is uniquely determined by its values on such standard
  intervals.
\end{remark}

\vspace{.2cm}

{\it Proof of Theorem~\ref{spectral subspaces}.} We let $\K$ denote the set
of sets of the form $B_1\times B_2$ with $B_1, B_2\in\B(\C)$.

Consider an arbitrary sequence of mutually disjoint sets from $\K$, $(B_1^{(k)}\times B_2^{(k)})_{k=1}^\infty$, and define
\begin{equation}\label{Br24}
  P_{S,T}\Big(\bigcup_{k=1}^\infty (B_1^{(k)}\times B_2^{(k)})\Big) :=
  \bigvee_{k=1}^\infty [P_S(B_1^{(k)})\wedge P_T(B_2^{(k)})].
\end{equation}
Clearly, $P:=P_{S,T}\Big(\bigcup_{k=1}^\infty (B_1^{(k)}\times
  B_2^{(k)})\Big)$ satisfies that
\begin{itemize}
  \item[(a)] $P$ is $S$- and $T$-invariant.
\end{itemize}
In addition, we prove that with $\CK = P(\CH)$,
\begin{itemize}
  \item[(b)] $\mu_{S|_{\CK},T|_{\CK}}$ is concentrated on
  $B:=\bigcup_{k=1}^\infty (B_1^{(k)}\times B_2^{(k)})$,
\end{itemize}
and
\begin{itemize}
  \item[(c)] $P$ is maximal w.r.t. the properties (a) and (b).
\end{itemize}
(c) will entail that the right-hand side of \eqref{Br24} is independent of
the way in which we write $B$ as a disjoint union of countably many
sets from $\K$, and hence, that $P_{S,T}(B)$ {\it does}, as indicated by the
notation, only depend on the set $B$.

To see that (b) holds, note that if $Q\in\CM$ is {\it any} $S$- and
$T$-invariant projection, and if we let $\CL=Q(\CH)$, then by \eqref{Br4} and
Lemma~\ref{intersect-inv},
\begin{eqnarray}\label{Br25}
\mu_{S|_{\CL},T|_{\CL}}(B)&=& \sum_{k=1}^\infty\mu_{S|_{\CL},T|_{\CL}}(B_1^{(k)}\times B_2^{(k)})\nonumber\\
&=& \sum_{k=1}^\infty \tau_{Q\CM Q}(P_{S|_{\CL}}(B_1^{(k)})\wedge
P_{T|_{\CL}}(B_2^{(k)}))\nonumber\\
&=& \sum_{k=1}^\infty \tau_{Q\CM Q}(P_S(B_1^{(k)})\wedge
P_T(B_2^{(k)})\wedge Q).
\end{eqnarray}

Then using Proposition~\ref{infinsum}, we get that
\begin{eqnarray*}
\mu_{S|_{\CK},T|_{\CK}}(B)&=& \sum_{k=1}^\infty \tau_{P\CM
  P}(P_S(B_1^{(k)})\wedge P_T(B_2^{(k)})\wedge P)\\
&=& \sum_{k=1}^\infty \tau_{P\CM P}(P_S(B_1^{(k)})\wedge
P_T(B_2^{(k)}))\\
&=& \frac{1}{\tau(P)}\sum_{k=1}^\infty \tr(e_S(B_1^{(k)})e_T(B_2^{(k)}))\\
&=& \frac{1}{\tau(P)}\tr\Big(\sum_{k=1}^\infty
e_S(B_1^{(k)})e_T(B_2^{(k)}) \Big)\\
&=&  \frac{1}{\tau(P)} \tau\Big(\bigvee_{k=1}^\infty [P_S(B_1^{(k)})\wedge
P_T(B_2^{(k)})]\Big)\\
&=& 1.
\end{eqnarray*}
Thus, (b) holds.

Now, suppose that $Q\in\CM$ is an $S$- and $T$-invariant projection, and
that $\mu_{S|_{\CL},T|_{\CL}}$ is concentrated on $B$, where
$\CL=Q(\CH)$. Then by Lemma~\ref{intersect-inv} and Proposition~\ref{infinsum},
\begin{eqnarray*}
  P\wedge Q &=& \Big(\bigvee_{k=1}^\infty [P_S(B_1^{(k)})\wedge
P_T(B_2^{(k)})]\Big)\wedge Q\\
&=& \bigvee_{k=1}^\infty [P_{S|_{\CL}}(B_1^{(k)})\wedge
P_{T|_{\CL}}(B_2^{(k)})]\\
&=& P_{\range\big(\sum_{k=1}^\infty e_{S|_{\CL}}(B_1^{(k)})
e_{T|_{\CL}}(B_2^{(k)})\big)}.
\end{eqnarray*}
Hence, Proposition~\ref{infinsum} and \eqref{Br25} imply that
\begin{eqnarray*}
  \tau_{Q\CM Q}(P\wedge Q)&=& \tr_{Q\CM Q}\Big(\sum_{k=1}^\infty
  e_{S|_{\CL}}(B_1^{(k)}) e_{T|_{\CL}}(B_2^{(k)})\Big)\\
&=& \sum_{k=1}^\infty \tr_{Q\CM Q}(e_{S|_{\CL}}(B_1^{(k)})
  e_{T|_{\CL}}(B_2^{(k)}))\\
  &=& \sum_{k=1}^\infty \tau_{Q\CM Q}(P_{S|_{\CL}}(B_1^{(k)})\wedge
P_{T|_{\CL}}(B_2^{(k)}))\\
&=& \mu_{S|_{\CL},T|_{\CL} }(B)\\
&=& 1.
\end{eqnarray*}
Thus, $P\wedge Q=Q$, and this shows that (c) holds.

As mentioned in Remark~\ref{open sets}, every open set $U\subseteq\C^2$ may
be written as a union of countably many mutually disjoint sets from
$\K$. Thus, we have now proved existence of $P_{S,T}(U)$ for every such
$U$, and for general $B\in \B(\C^2)$ we will define
\begin{equation}\label{Br26}
  P_{S,T}(B):= \bigwedge_{B\subseteq U,\; U\subseteq \C^2\;{\rm open}}
  P_{S,T}(U).
\end{equation}
Then again, $P:=P_{S,T}(B)$ satisfies that 
\begin{itemize}
  \item[(a)] $P$ is $S$- and $T$-invariant.
\end{itemize}
Moreover, we prove that with $\CK = P(\CH)$,
\begin{itemize}
  \item[(b)] $\mu_{S|_{\CK},T|_{\CK}}$ is concentrated on
  $B$,
\end{itemize}
and
\begin{itemize}
  \item[(c)] $P$ is maximal w.r.t. the properties (a) and (b).
\end{itemize}
These propterties will entail that when $B$ happens to be a union of countably many mutually disjoint sets from
$\K$, then \eqref{Br26} agrees with the previous definition of
$P_{S,T}(B)$ (cf. \eqref{Br24}).

Now, to see that (b) holds, note that $\mu_{S|_{\CK},T|_{\CK}}$ is regular
(cf. \cite[Theorem~7.8]{Fo}), and hence
\begin{equation}\label{Br27}
\mu_{S|_{\CK},T|_{\CK}}(B) =
\inf\{\mu_{S|_{\CK},T|_{\CK}}(U)\,|\,B\subseteq U,\; U\subseteq \C^2\;{\rm
  open}\}.
\end{equation}
Let $U$ be any open subset of $\C^2$ containing $B$. Write $U$ as a union
of countably many mutually disjoint sets from $\K$:
\[
U=\bigcup_{k=1}^\infty (B_1^{(k)}\times B_2^{(k)}).
\]
Then, according to \eqref{Br25},
\begin{eqnarray*}
\mu_{S|_{\CK},T|_{\CK}}(U)
&=& \sum_{k=1}^\infty \tau_{P\CM P}(P_S(B_1^{(k)})\wedge
P_T(B_2^{(k)})\wedge P),
\end{eqnarray*}
and using Proposition~\ref{infinsum} and Lemma~\ref{intersect-inv} we find that
\begin{eqnarray*}
\mu_{S|_{\CK},T|_{\CK}}(U)&=& \tr_{P\CM P}\Big(\sum_{k=1}^\infty
  e_{S|_{\CK}}(B_1^{(k)}) e_{T|_{\CK}}(B_2^{(k)})\Big)\\
  &=& \tau_{P\CM P}\Big(\bigvee_{k=1}^\infty [P_{S|_{\CK}}(B_1^{(k)})\wedge
  P_{T|_{\CK}}(B_2^{(k)})]\Big)\\
  &=& \tau_{P\CM P}(P_{S,T}(U)\wedge P)\\
  &=& \tau_{P\CM P}(P)\\
  &=& 1,
\end{eqnarray*}
where $P_{S,T}(U)$ is given by \eqref{Br24}. Hence by
\eqref{Br27}, $\mu_{S|_{\CK},T|_{\CK}}$ is concentrated on $B$.

Finally, if $Q\in\CM$ is any $S$- and $T$-invariant projection, and if
$\mu_{S|_{\CL},T|_{\CL}}$ is concentrated on $B$, where $\CL=Q(\CH)$, then
$\mu_{S|_{\CL},T|_{\CL}}$ is concentrated on $U$ for every open set $U$
containing $B$. Hence, by the first part of the proof, $Q\leq P_{S,T}(U)$
for every such $U$, and it follows from the definition of $P_{S,T}(B)$ that
$Q\leq P$.

Concerning \eqref{measure-trace}, note that if $B=B_1\times B_2$, where
$B_1, B_2 \in\B(\C)$, then, by the definitions of $\mu_{S,T}$ and $P_{S,T}(B)$,
\eqref{measure-trace} holds. If $B$ is a disjoint union of sets $
    (B^{(k)})_{k=1}^\infty=(B_1^{(k)}\times B_2^{(k)})_{k=1}^\infty$, where
    $B_i^{(k)}\in\B(\C)$, $k\in\N$, $i=1, 2$, then
\begin{eqnarray*}
  \mu_{S,T}(B)&=& \sum_{k=1}^\infty \tau(P_{S,T}(B_1^{(k)}\times
  B_2^{(k)}))\\
  &=&  \sum_{k=1}^\infty \tau(\supp(e_S(B_1^{(k)})e_T(B_2^{(k)}))).
\end{eqnarray*}
Applying Proposition~\ref{infinsum} we thus find that
\begin{eqnarray*}
  \mu_{S,T}(B)&=& \tau\Big(\supp\Big( \sum_{k=1}^\infty
  e_S(B_1^{(k)})e_T(B_2^{(k)})\Big)\Big)\\
  &=& \tau\Big(\bigvee_{k=1}^\infty
  \supp(e_S(B_1^{(k)})e_T(B_2^{(k)}))\Big)\\
  &=& \tau\Big(\bigvee_{k=1}^\infty
  P_{S,T}(B_1^{(k)}\times B_2^{(k)}))\Big)\\
&=& \tau(P_{S,T}(B)).
\end{eqnarray*}

Finally, for general $B\in\B(\C^2)$, since $\mu_{S,T}$ is regular,
\begin{eqnarray*}
\mu_{S,T}(B)&= & \inf\{\mu_{S,T}(U)\,|\,B\subseteq U\subseteq \C^2,\,
U \,{\rm open}\}\\
&=& \inf\{\tau(P_{S,T}(U))\,|\,B\subseteq U\subseteq \C^2,\,
U \,{\rm open}\}\\
&=& \tau\Big(\bigwedge_{B\subseteq U\subseteq \C^2,\,
U \,{\rm open}}P_{S,T}(U)\Big)\\
&=& \tau(P_{S,T}(B)). \qquad\qquad \endproof
\end{eqnarray*}

\vspace{.2cm}

The proof given above may clearly be generalized to the case of $n$ commuting operators $T_1, \ldots, T_n\in\CM$, so
  that Theorem~\ref{spectral subspaces} has a slightly more general
  version:

\begin{thm}\label{spectral subspaces n} Let $n\in\N$, let $T_1,\ldots, T_n\in\CM$ be commuting
  operators, and let $B\subseteq \C^n$ be any Borel set. Then there is a
  maximal closed subspace, $\CK = \CK_{T_1,\ldots, T_n}(B)$, affiliated
  with $\CM$ which is $T_i$-invariant for every $i\in\{1,
  \ldots, n\}$, and such  
  that the Brown measure $\mu_{T_1|_\CK,\ldots, T_n|_\CK}$ is concentrated on
  $B$. Let $P_{T_1,\ldots, T_n}(B)\in\CM$ denote the projection onto
  $\CK_{T_1,\ldots, T_n}(B)$. Then more precisely,
  \begin{itemize}
    \item[(i)] if $B=B_1\times\cdots\times B_n$ with $B_i\in\B(\C)$, then
  \begin{equation}
    P_{T_1,\ldots, T_n}(B)=\bigwedge_{i=1}^n P_{T_i}(B_i),
  \end{equation}
    \item[(ii)] if $B$ is a disjoint union of sets $
    (B^{(k)})_{k=1}^\infty=(B_1^{(k)}\times\cdots\times B_n^{(k)})_{k=1}^\infty$, where
    $B_i^{(k)}\in\B(\C)$, $k\in\N$, $i=1,\ldots, n$, then
    \begin{equation}
      P_{T_1,\ldots, T_n}(B)=\bigvee_{k=1}^\infty P_{T_1,\ldots, T_n,}(B^{(k)}),
    \end{equation}
    \item[(iii)] and for general $B\in\B(\C^n)$,
      \begin{equation}
      P_{T_1,\ldots, T_n}(B)=\bigwedge_{B\subseteq U,\;U\subseteq \C^n \;{\rm open}}P_{T_1,\ldots, T_n}(U).
      \end{equation}
  \end{itemize}
  Moreover, for every  $B\in\B(\C^n)$,
  \begin{equation}
    \mu_{T_1,\ldots, T_n}(B)=\tau(P_{T_1,\ldots, T_n}(B)).
  \end{equation}
\end{thm}

\section{An alternative characterization of $\mu_{S,T}$.}

In this final section we are going to give a characterization of the
Brown measure of two commuting operators in $\CM$, which is
different from the one we gave in Theorem~\ref{commuting}. Recall from
\cite{Bro} that for $T\in\CM$, the Brown measure of $T$, $\mu_T$, is the unique
compactly supported Borel probability measure on $\C$ which satisfies
the identity
\begin{equation}
\tau(\log|T-\lambda\unit|)=\int_\C \log|z-\lambda|\,\d\mu_T(z)
\end{equation}
for all $\lambda\in\C$. 

We are going to prove that a similar property characterizes $\mu_{S,T}$:

\begin{thm}\label{uniqueness} Let $S,T\in\CM$ be commuting operators. Then
  $\mu_{S,T}$ is the unique compactly supported Borel probability measure on
  $\C^2$ which satisfies the identity
  \begin{equation}\label{alphaS+betaT}
    \tau(\log|\alpha S + \beta T -\unit|) = \int_{\C^2}\log|\alpha z
    + \beta w - 1|\,\d\mu_{S,T}(z,w).
  \end{equation}
for all $\alpha, \beta \in\C$.
\end{thm}

\vspace{.2cm}

\begin{remark}\label{equivalent} Let $S,T\in\CM$ be as in
  Theorem~\ref{uniqueness}. Note that if $\mu_{S,T}$ satisfies
  \eqref{alphaS+betaT} for all $\alpha, \beta \in\C$, then for all $\alpha,
  \beta,\lambda \in\C$,
   \begin{equation}\label{corr1}
    \tau(\log|\alpha S + \beta T -\lambda\unit|) = \int_{\C^2}\log|\alpha z
    + \beta w - \lambda|\,\d\mu_{S,T}(z,w).
  \end{equation}
  This is clear for $\lambda\neq 0$, and for $\lambda=0$, \eqref{corr1} follows from the fact that two subharmonic functions defined in $\C$ coincide iff they agree almost everywhere w.r.t. Lebesgue measure. It nowfollows from Brown's characterization of $\mu_{\alpha S+\beta T}$
  that $\mu_{\alpha S+\beta T}$ is the push-forward measure
  $\nu_{\alpha,\beta}$ of $\mu_{S,T}$ via the map $(z,w)\mapsto \alpha z +
  \beta w$. On the other hand, \emph{if}  $\nu_{\alpha, \beta} = \mu_{\alpha S +
  \beta T}$, then \eqref{alphaS+betaT} holds. 
\end{remark}

\vspace{.2cm}

Recall from \cite{HS} that
the {\it modified spectral radius} of $T\in\CM$, $r'(T)$, is defined by
\begin{equation}
r'(T):= \max\{|z|\,|\, z\in \supp(\mu_T)\}.
\end{equation}
Also recall from \cite[Corollary~2.6]{HS} that in fact
\begin{equation}\label{r'(T)}
  r'(T)=\lim_{p\rightarrow
  \infty}\Big(\lim_{n\rightarrow\infty}\|T^n\|_{\frac pn}^{\frac1n}\Big).
\end{equation}

\begin{lemma}\label{spectralradii} Let $S,T\in\CM$ be commuting
  operators. Then the modified spectral radii, $r'(S)$, $r'(T)$, $r'(ST)$ and
  $r'(S+T)$, satisfy the inequalities
  \begin{equation}\label{times}
  r'(ST)\leq r'(S)\cdot r'(T),
  \end{equation}
  and
  \begin{equation}\label{plus}
  r'(S+T)\leq r'(S)+r'(T).
  \end{equation}
\end{lemma}
  
\proof \eqref{times} follows from \eqref{r'(T)} and the generalized H\"
older inequality (cf. \cite{FK}): For $A,B\in\CM$ and for $0<p,q,r\leq
\infty$ with $\frac1r =\frac1p + \frac 1q$,
\[
\|AB\|_r\leq \|A\|_p\|B\|_q.
\]

To prove  \eqref{plus}, note that
\[
\supp(\mu_S)\subseteq \{z\in\C\,|\,\re z\leq r'(S)\},
\]
and
\[
\supp(\mu_T)\subseteq \{z\in\C\,|\,\re z\leq r'(T)\}.
\]
According to \cite[Theorem~4.1]{Bro}, $\mu_{\e^S}$ and $\mu_{\e^T}$ are the
push-forward measures of $\mu_S$ and $\mu_T$, respectively, via the
map $z\mapsto \e^z$. Hence,
\[
\supp(\mu_{\e^S})\subseteq \{z\in\C\,|\,|z|\leq \e^{r'(S)}\},
\]
and
\[
\supp(\mu_{\e^T})\subseteq \{z\in\C\,|\,|z|\leq \e^{r'(T)}\},
\]
and it follows from \eqref{times} that
\begin{eqnarray*}
r'(\e^{S+T})&=& r'(\e^S\e^T)\\
&\leq &\e^{r'(S)}\e^{r'(T)} \\
&=& \e^{r'(S)+r'(T)}.
\end{eqnarray*}
Thus, $\supp(\mu_{\e^{S+T}})\subseteq \overline{B(0,\e^{r'(S)+r'(T)})}$,
and then, by one more application of \cite[Theorem~4.1]{Bro},
\[
\supp(\mu_{S+T})\subseteq  \{z\in\C\,|\,\re z\leq r'(S)+ r'(T)\}.
\]
Repeating this argument, we find that for arbitrary $\theta\in [0, 2\pi[$,
\[
\supp(\mu_{\e^{\i\theta}(S+T)})\subseteq  \{z\in\C\,|\,\re z\leq r'(S)+ r'(T)\},
\]
i.e.
\[
\supp(\mu_{S+T})\subseteq  \{z\in\C\,|\,\re (\e^{-\i\theta}z)\leq r'(S)+
r'(T)\}.
\]
Since $\theta$ was arbitrary, we conclude that
\[
\supp(\mu_{S+T})\subseteq \overline{B(0, r'(S)+r'(T))},
\]
and this proves \eqref{plus}. $\endproof$

\vspace{.2cm}

\begin{lemma}\label{scalarmult}Let $S,T\in\CM$ be commuting
  operators, and let $\alpha, \beta\in\C$. Then $\mu_{\alpha
    S, \beta T}$ is the push-forward measure of $\mu_{S,T}$ via the map
    $h_{\alpha,\beta}: \C\times \C\rightarrow \C\times \C$ given by
    \[
    h_{\alpha,\beta}(z,w)= (\alpha z , \beta w).
    \]
\end{lemma}

\proof Recall that $\mu_{\alpha S, \beta T}$ is uniquely determined by the
property that for all $B_1, B_2\in \B(\C)$,
\begin{equation}\label{eq1}
\mu_{\alpha S, \beta T}(B_1\times B_2) = \tau(P_{\alpha S}(B_1)\wedge
P_{\beta T}(B_2)).
\end{equation}
Now, it is easily seen that for $\alpha\neq 0$ and $\beta\neq 0$,
$P_{\alpha S}(B_1)= P_S(\frac1\alpha B_1)$ and $P_{\beta T}(B_2)=
P_T(\frac1\beta B_2)$. Hence,
\begin{eqnarray}
\mu_{\alpha S, \beta T}(B_1\times B_2) &=& \tau(P_{S}(\textstyle{\frac1\alpha}B_1)\wedge
P_{T}(\textstyle{\frac1\beta}B_2))\nonumber \\
&=& \mu_{S,T}(\textstyle{\frac1\alpha}B_1\times \textstyle{\frac1\beta}B_2
)\nonumber \\
&=& \mu_{S,T}(h_{\alpha,\beta}^{-1}(B_1\times B_2)). \label{eq2}
\end{eqnarray}
If for instance $\alpha =0$, then $P_{\alpha S}(B_1)=0$ if $0\notin B_1$
and  $P_{\alpha S}(B_1)=\unit$ if $0\in B_1$. It then follows 
that \eqref{eq2} holds in this case as well. Similar arguments apply if $\beta =0$. $\endproof$

\vspace{.2cm}

{\it Proof of Theorem~\ref{uniqueness}.} As noted in
Remark~\ref{equivalent}, it suffices to prove that for all $\alpha,
\beta\in\C$, $\mu_{\alpha S + \beta T}$ is the push-forward measure of
$\mu_{S,T}$ via the map $(z,w)\mapsto \alpha z + \beta w$. At first we
will consider the case $\alpha=\beta =1$. Define $a: \C\times \C\rightarrow
\C$ by
\[
a(z,w)=z+w, \quad (z,w\in\C).
\]
We are going to prove that for all $B\in\B(\C)$,
\begin{equation}\label{push-forward}
\mu_{S+T}(B)=\mu_{S,T}(a^{-1}(B)).
\end{equation}
It suffices to show that for every open set $U\subseteq \C$, 
\begin{equation}
\mu_{S+T}(U)\geq \mu_{S,T}(a^{-1}(U)).
\end{equation}
Indeed, if this holds, then by regularity of $\mu_{S+T}$ and $a(\mu_{S,T})$, for every Borel set $B\subseteq \C$,
\begin{eqnarray*}
\mu_{S+T}(B)&=& \inf\{\mu_{S+T}(U)\,|\,B\subseteq U,\,U\,{\rm open}\}\\
&\geq & \inf\{a(\mu_{S,T})(U)\,|\,B\subseteq U,\,U\,{\rm open}\}\\
&=& \mu_{S,T}(a^{-1}B).
\end{eqnarray*}
Since both measures are probability measures, and the above inequality holds for both $B$ and $B^c$, we must have identity. That is, \eqref{push-forward} holds. 

Now, let $U\subseteq \C$ be any open set. Then $V:=a^{-1}(U)$ is open in $\C^2$ and we may write $V$ as a countable union of mutually disjoint "boxes",
\[
V= \bigcup_{n=1}^\infty I(z_n,\delta_n)\times I(w_n,\delta_n),
\]
where for $z\in\C$ and $\delta>0$,
\begin{equation}\label{p1}
I(z,\delta):= \{w\in\C\,|\, \re(z)-\delta < \re(w)\leq \re(z)+\delta,\; \im(z)-\delta < \im(w)\leq \im(z)+\delta \}.
\end{equation}
We can even choose $z_n, w_n$ and $\delta_n>0$ so that
\begin{equation}\label{ball-inclusion}
\overline{B(z_n+w_n, 2\sqrt2\delta_n)}\subseteq U.
\end{equation}
This requires a little consideration and for the convenience of the reader, we provide an argument in Lemma~\ref{box-union} below.

Now, according to Theorem~\ref{spectral subspaces},
\begin{eqnarray*}
\mu_{S,T}(a^{-1}(U))&=&\tau(P_{S,T}(V))\\
&=& \tau\Big(\bigvee_{n=1}^\infty \big[P_S(I(z_n,\delta_n))\wedge P_T(I(w_n,\delta_n))\big]\Big),
\end{eqnarray*}
and we also have that
\[
\mu_{S+T}(U)=\tau(P_{S+T}(U)). 
\]
Hence, it suffices to prove that for every $n\in\N$,
\begin{equation}
P_{S+T}(U)\geq P_S(I(z_n,\delta_n))\wedge P_T(I(w_n,\delta_n)).
\end{equation}
Fix $n\in\N$, and set $P=P_S(I(z_n,\delta_n))\wedge P_T(I(w_n,\delta_n))$. Then by Lemma~\ref{spectralradii},
\begin{eqnarray*}
r'\big([S+T-(z_n+w_n)\unit]|_{P(\CH)}\big)&\leq& r'\big([S-z_n\unit]|_{P(\CH)}\big)+r'\big([T-w_n\unit]|_{P(\CH)}\big)\\
&\leq& r'\big([S-z_n\unit]|_{P_S(I(z_n,\delta_n))(\CH)}\big)+r'\big([T-w_n\unit]|_{P_T(I(w_n,\delta_n))(\CH)}\big)\\
&\leq & 2\sqrt2 \delta_n,
\end{eqnarray*}
and it follows that $\mu_{S+T|_{P(\CH)}}$ is concentrated on $\overline{B(z_n+w_n, 2\sqrt2 \delta_n)}\subseteq U$. Hence, $P\leq P_{S+T}(U)$, and we are done.

Now, if $\alpha,\beta\in\C$, then we conclude from the above and
Lemma~\ref{scalarmult} that
\begin{eqnarray*}
\mu_{\alpha S + \beta T}(B)&=& \mu_{\alpha S, \beta T}(a^{-1}(B))\\
 &=& \mu_{S,T}(h_{\alpha,\beta}^{-1}(a^{-1}(B)))\\
 &=& \mu_{S,T}((a\circ h_{\alpha,\beta})^{-1}(B)),
\end{eqnarray*}
and since $(a\circ h_{\alpha,\beta})(z,w)=\alpha z + \beta w$, this completes
the proof of the identity \eqref{alphaS+betaT}.

To prove uniqueness of $\mu_{S,T}$, suppose that $\nu$ is a compactly
supported Borel probability measure on $\C^2$ which satisfies the
identity \eqref{alphaS+betaT} for all $\alpha, \beta\in\C$. That is,  for all $\alpha,
\beta\in\C$, $\mu_{\alpha S + \beta T}$ is the push-forward measure of
$\nu$ via the map $(z,w)\mapsto \alpha z + \beta w$. Then, to prove that
$\nu = \mu_{S,T}$, it suffices to prove that for all $y= (y_1, \ldots,
y_4)\in \R^4$,
\begin{equation}
  \int_{\R^4}\e^{\i(y,x)}\,\d\mu_{S,T}(x)=  \int_{\R^4}\e^{\i(y,x)}\,\d\nu(x)
\end{equation}
(here we identify $\C$ with $\R^2$). For $x= (x_1, \ldots,
x_4)\in \R^4$ and $y= (y_1, \ldots, y_4)\in \R^4$, note that
\[
(y,x)= \re\Big((y_1-\i y_2)(x_1+\i x_2)+(y_3-\i y_4)(x_3+\i x_4)\Big),
\]
and hence with $\alpha = y_1-\i y_2$ and $\beta = y_3-\i y_4$ we find that
\begin{eqnarray*}
  \int_{\R^4}\e^{\i(y,x)}\,\d\mu_{S,T}(x)&=&  \int_{\C^2}\e^{\i \re(\alpha
  z+ \beta w)}\,\d\mu_{S,T}(z,w)\\
&=& \int_{\C}\e^{\i \re z}\,\d\mu_{\alpha S + \beta T}(z)\\
&=& \int_{\C^2}\e^{\i \re(\alpha
  z+ \beta w)}\,\d\nu(z,w)\\\\
&=& \int_{\R^4}\e^{\i(y,x)}\,\d\nu(x),
\end{eqnarray*}
as desired. $\endproof$

\vspace{.2cm}

\begin{remark} In the proof above it was shown that for $U\subseteq \C$ an open set, we have the following inequality:
\begin{equation}\label{sept18}
P_{S+T}(U)\geq P_{S,T}(a^{-1}(U)).
\end{equation}
But it was also shown that the two projections above have the same trace:
\[
\tau(P_{S+T}(U))=\mu_{S+T}(U)=\mu_{S,T}(a^{-1}(U))=\tau(P_{S,T}(a^{-1}(U))).
\]
Hence, the two projections in \eqref{sept18} are identical, and by Theorem~\ref{spectral subspaces}~(iii), for \emph{every} Borel set $B\subseteq \C$, we must have that
\begin{equation}\label{add3}
P_{S+T}(B)= P_{S,T}(p^{-1}(B)).
\end{equation}
\end{remark}

\vspace{.2cm}

As in the previous section, one can easily generalize the proof given above
to the case of an arbitrary finite set of commuting operators, $\{T_1,
  \ldots, T_n\}$. That is, we actually have the following alternative
  description of $\mu_{T_1,\ldots, T_n}$:

\begin{thm}\label{n commuting operators} Let $n\in\N$, and let $T_1, \ldots, T_n$ be mutually commuting
  operators in $\CM$. Then
  $\mu_{T_1, \ldots, T_n}$ is the unique compactly supported Borel probability measure on
  $\C^n$ which satisfies the identity
  \begin{equation}
    \tau(\log|\alpha_1 T_1 + \ldots + \alpha_n T_n -\unit|) =
    \int_{\C^n}\log|\alpha_1 z_1 + \ldots + \alpha_n z_n - 1|\,\d\mu_{T_1,
  \ldots, T_n}(z_1,\ldots, z_n)
  \end{equation}
for all $\alpha_1,\ldots, \alpha_n \in\C$.
\end{thm}

\vspace{.2cm}

\begin{lemma}\label{box-union} Define $a:\C^2\rightarrow\C$ by 
\[
a(z,w)=z+w,
\]
and let $U\subseteq \C$ be an open set. Then for every pair $(S,T)$ of commuting operators in $\CM$, we may write $V:=a^{-1}(U)$ as a countable disjoint union of sets $\big(I(z_n,\delta_n)\times I(w_n,\delta_n)\big)_{n=1}^\infty$, where for $z\in\C$ and $\delta>0$,
\[
I(z,\delta):= \{w\in\C\,|\, \re(z)-\delta < \re(w)\leq \re(z)+\delta,\; \im(z)-\delta < \im(w)\leq \im(z)+\delta \}.
\]
Moreover, we can ensure that for each $n\in\N$,
\begin{equation}\label{conditionp}
\overline{B(z_n+w_n, 2\sqrt2\delta_n)}\subseteq U.
\end{equation}
\end{lemma}

\proof Divide $\C^2$ into mutually disjoint "boxes" of the form $I(z,1)\times I(w,1)$. Take the countably many of these, $\big(I(z_n^{(1)},1)\times I(w_n^{(1)},1)\big)$, which are contained in $V$ and satisfy \eqref{conditionp} (with $\delta_n=1$). Next, divide $\C^2$ into boxes once more by cutting each of the previous ones into 16 boxes of equal size (edge lenght $\frac12$). Take those, $\big(I(z_n^{(2)},\frac12)\times I(w_n^{(2)},\frac12)\big)$, which are not contained in $\bigcup_n I(z_n^{(1)},1)\times I(w_n^{(1)},1)$ but which are contained in $V$ and satisfy \eqref{conditionp} (with $\delta_n=\frac12$). Continue like this and obtain a set $V_0\subseteq V$ of the form
\[
V_0 = \bigcup_{m=1}^\infty \bigcup_n I(z_n^{(m)},2^{-m})\times I(w_n^{(m)},2^ {-m}).
\]
We claim that this disjoint union is all of $V$. Indeed, let $(z,w)\in V$. Then $z+w\in U$ so there is an $\eps>0$ such that $\overline{B(z+w,\eps)}\subseteq U$. Choose $m\in\N$ so large that $2^{-m}< \frac{\eps}{4\sqrt2}$. Then, if $|z-z'|\leq \sqrt2 2^{-m}$ and $|w-w'|\leq \sqrt2 2^{-m}$, one has that
\[
\overline{B(z'+w',2\sqrt2 2^{-m}})\subseteq \overline{B(z+w, 4\sqrt 2 2^{-m})}\subseteq U.
\]
Thus, when we have divided $\C^2$ into boxes of edge length at most $2^{-m}$, then the one box containing $(z,w)$ will satisfy \eqref{conditionp}. Then it is just a matter of taking $m$ so large that the box is also contained in $V$. It follows that $(z,w)\in V_0$. $\endproof$

\vspace{.2cm}

\begin{remark} Let $S\in \CM$ be invertible, and let $T_1, \ldots, T_n$ be
  mutually commuting operators in $\CM$. Then
  \begin{equation}
    \mu_{S^{-1}T_1S, \ldots,S^{-1}T_nS}= \mu_{T_1,\ldots, T_n}.
  \end{equation}
Indeed, this follows from the characterization of $\mu_{T_1,\ldots, T_n}$
given in Theorem~\ref{n commuting operators} and from the fact that for all
$T\in \CM$, $\mu_{S^{-1}TS}= \mu_T$ (cf. \cite{Bro}).
\end{remark}

\vspace{.2cm}

\begin{prop}\label{product} Let $S,T\in\CM$ be commuting
  operators. Then $\mu_{ST}$ is the push-forward measure of $\mu_{S,T}$ via
  the map $m: (z,w)\mapsto zw$.
\end{prop}

\proof The proof is essentially the same as the one we gave above when considering the map $a: (z,w)\mapsto z+w$. Again it suffices to show that for every open set $U\subseteq \C$, 
\begin{equation}
\mu_{ST}(U)\geq \mu_{S,T}(m^{-1}(U)),
\end{equation}
and for such an open set $U$ we write $V:=m^{-1}(U)$ as a countable union of mutually disjoint "boxes" as in \eqref{p1}, but this time we make sure that $\delta_n>0$ is so small that 
\begin{equation}
\overline{B(z_nw_n, \sqrt2\delta_n(\|T\|+|z_n|))}\subseteq U.
\end{equation}
As in the previous case, one only has to show that for every $n\in\N$,
\[
P_{ST}(U)\geq P_S(I(z_n,\delta_n))\wedge P_T(I(w_n,\delta_n)).
\]
Fix $n\in\N$ and set $P=P_S(I(z_n,\delta_n))\wedge P_T(I(w_n,\delta_n))$. Since
\[
r'\big([S-z_n\unit]|_{P(\CH)}\big)\leq r'\big([S-z_n\unit]|_{P_S(I(z_n,\delta_n))(\CH)}\big)\leq \sqrt2\delta_n,
\]
and
\[
r'\big([T-w_n\unit]|_{P(\CH)}\big)\leq r'\big([T-w_n\unit]|_{P_T(I(w_n,\delta_n))(\CH)}\big)\leq \sqrt2\delta_n,
\]
and since
\[
ST-z_nw_n\unit= (S-z_n\unit)T + z_n(T-w_n\unit),
\]
we have (cf. Lemma~\ref{spectralradii}) that
\begin{eqnarray*}
r'\big([ST-z_nw_n\unit]|_{P(\CH)})&\leq & r'\big([(S-z_n\unit)T]|_{P(\CH)}\big)+ |z_n|r'\big([T-w_n\unit]|_{P(\CH)}\big)\\
&\leq & r'\big([S-z_n\unit]|_{P(\CH)}\big)\|T\|+|z_n|r'\big([T-w_n\unit]|_{P(\CH)}\big)\\
&\leq & \sqrt2\delta_n(\|T\|+|z_n|).
\end{eqnarray*}
Thus, $\mu_{ST|_{P(\CH)}}$ is concentrated on $\overline{B(z_nw_n, \sqrt2\delta_n(\|T\|+|z_n|))}\subseteq U$, and therefore $P\leq P_{ST}(U)$, as desired. $\endproof$

\vspace{.2cm}

\begin{remark} As in the additive case, we infer from the proof given above that for every Borel set $B\subseteq \C$ we have:
\begin{equation}\label{pol6}
P_{ST}(B)=P_{S,T}(m^{-1}(B)).
\end{equation}
\end{remark}

\vspace{.2cm}

\begin{prop} Consider type II$_1$ factors $\CM_1$ and $\CM_2$ with faithful
  tracial states $\tau_1$ and $\tau_2$, respectively. Let $S\in\CM_1$ and
  $T\in\CM_2$. Then
  \begin{equation}\label{tensor}
    \mu_{S\otimes\unit, \unit\otimes T} = \mu_S\otimes \mu_T,
\end{equation}
and it follows that
  \begin{eqnarray}\label{convolution}
    \mu_{S\otimes\unit + \unit\otimes T}&=& \mu_S \ast \mu_T,\label{p2}\\
    \mu_{S\otimes T} &=& \mu_S \star \mu_T,\label{m1}
  \end{eqnarray}
  where $\ast$ ($\star$, resp.) denotes additive (multiplicative, resp.) convolution.
\end{prop}

\proof $\mu_{S\otimes\unit + \unit\otimes T}$ ($\mu_{S\otimes T}$, resp.) is the push-forward measure
of $\mu_{S\otimes\unit, \unit\otimes T}$ via the map $a: (z,w)\mapsto z+w$ ($m: (z,w)\mapsto zw)$, resp.), and $\mu_S \ast \mu_T$ ($\mu_S \star \mu_T$, resp.) is the push-forward measure
of $\mu_S\otimes \mu_T$ via that same map. Thus, \eqref{p2} and \eqref{m1}
follow from \eqref{tensor}. To see that the latter holds, let $B_1,
B_2\in\B(\C)$. It is easily seen that
\[
P_{S\otimes\unit}(B_1)= P_S(B_1)\otimes\unit \quad {\rm and}
\quad P_{\unit\otimes T}(B_2)= \unit\otimes P_T(B_2).
\]
Hence,
\begin{eqnarray*}
  \mu_{S\otimes\unit, \unit\otimes T}(B_1\times B_2)&=& (\tau_1\otimes
  \tau_2)[(P_S(B_1)\otimes\unit)\cap (\unit\otimes P_T(B_2))]\\
  &=& \tau_1(P_S(B_1))\tau_2(P_T(B_2))\\
  &=& \mu_S(B_1)\mu_T(B_2)\\
  &=& (\mu_S\otimes \mu_T)(B_1\times B_2).
\end{eqnarray*}
This proves \eqref{tensor}. $\endproof$

\section{Polynomials in $n$ commuting variables.}

In this final section we will prove:

\begin{thm}\label{polynomial} 
Let $n\in\N$, and let $q$ be a polynomial in $n$ commuting variables, i.e. $q\in\C[z_1,\ldots,z_n]$. Then for every $n$-tuple $(T_1,\ldots,T_n)$ of commuting operators in $\CM$, one has that
\begin{equation}
\mu_{q(T_1,\ldots,T_n)}=q(\mu_{T_1,\ldots,T_n}),
\end{equation}
where $q(\mu_{T_1,\ldots,T_n})$ is the push-forward measure of $\mu_{T_1,\ldots,T_n}$ via $q:\C^n\rightarrow\C$.
\end{thm}

The proof relies on the previous sections and a few technical lemmas. 

\begin{lemma}\label{pol1} 
Given $n\in\N$ and commuting operators $T_1,\ldots, T_n\in\CM$. Let $1\leq i < n$. Then for all Borel sets $A\subseteq\C^i$ and $B\subseteq \C^{n-i}$, one has that  
\begin{equation}\label{pol2}
P_{T_1,\ldots,T_n}(A\times B)= P_{T_1,\ldots,T_i}(A)\wedge P_{T_{i+1},\ldots,T_n}(B).
\end{equation}
\end{lemma} 

\proof Fix Borel sets $A_1,\ldots,A_i\subseteq \C$, put $A=A_1\times\cdots\times A_i$, and let $B\subseteq \C^{n-i}$ be any open set. Then we may write $B$ as a disjoint union of cartesian products of Borel sets $B_m^{(j)}$, $i+1\leq j\leq n$, i.e.
\[
B= \bigcup_{m=1}^\infty\big[B_m^{(i+1)}\times\cdots\times B_m^{(n)}\big].
\]
Then by Theorem~\ref{spectral subspaces n}, \eqref{pol2} holds. Now, let $B\subseteq \C^{n-i}$ be any Borel set. Then by Theorem~\ref{spectral subspaces n}, 
\[
P_{T_{i+1},\ldots,T_n}(B)=\bigwedge_{B\subseteq U,\;U\;{\rm open}}P_{T_{i+1},\ldots,T_n}(U),
\]
and it follows from the above that 
\begin{eqnarray}
P_{T_1,\ldots,T_i}(A)\wedge P_{T_{i+1},\ldots,T_n}(B)&=&\bigwedge_{B\subseteq U,\;U\;{\rm open}}P_{T_1,\ldots,T_i}(A)\wedge P_{T_{i+1},\ldots,T_n}(U)\\
&=& \bigwedge_{B\subseteq U,\;U\;{\rm open}}P_{T_1,\ldots,T_n}(A\times U)\label{pol3}.
\end{eqnarray}
Since $P_{T_1,\ldots,T_n}(A\times B)\leq P_{T_1,\ldots,T_n}(A\times U)$ when $B\subseteq U$, this implies that
\begin{equation}\label{pol4}
P_{T_1,\ldots,T_i}(A)\wedge P_{T_{i+1},\ldots,T_n}(B)\geq P_{T_1,\ldots,T_n}(A\times B).
\end{equation}
On the other hand, \eqref{pol3} shows that
\begin{eqnarray*}
\tau(P_{T_1,\ldots,T_i}(A)\wedge P_{T_{i+1},\ldots,T_n}(B))&=&\inf_{B\subseteq U,\;U\;{\rm open}}\tau(P_{T_1,\ldots,T_n}(A\times U))\\
&=&\inf_{B\subseteq U,\;U\;{\rm open}}\mu_{T_1,\ldots,T_n}(A\times U).
\end{eqnarray*}
For fixed $A$ as above, the map $B\mapsto \mu_{T_1,\ldots,T_n}(A\times B)$ defined on Borel subsets of $\C^{n-i}$ is a finite (hence regular) Borel measure. Therefore we conclude that
\begin{eqnarray*}
\tau(P_{T_1,\ldots,T_i}(A)\wedge P_{T_{i+1},\ldots,T_n}(B))&=& \mu_{T_1,\ldots,T_n}(A\times B)\\
&=& \tau(P_{T_1,\ldots,T_n}(A\times B)).
\end{eqnarray*}
Comparing this identity with \eqref{pol4}, we find that \eqref{pol2} holds for this particular $A$. Now fix an arbitrary Borel set $B\subseteq \C^{n-i}$, and proceed in the same manner. That is, at first assume that $A\subseteq \C^i$ is open and verify that \eqref{pol2} holds in this case. Then finally consider an arbitrary Borel set $A$. $\endproof$  

\vspace{.2cm}

\begin{lemma}\label{pol5} Let $n\in\N$ and let $\alpha\in\C$. Define $a_n,\,m_n^{(\alpha)}:\C^{n+1}\rightarrow\C$ by 
\begin{equation}
a_n(z_1,\ldots,z_n,z_{n+1})=(z_1,\ldots,z_n+z_{n+1}),
\end{equation}
\begin{equation}
m_n^{(\alpha)}(z_1,\ldots,z_n,z_{n+1})=(z_1,\ldots,\alpha z_nz_{n+1}).
\end{equation}
Then for any $(n+1)$-tuple $(T_1,\ldots,T_{n+1})$ of commuting operators in $\CM$ one has that
\begin{equation}\label{pol7}
\mu_{T_1,\ldots, T_{n-1},T_n+T_{n+1}}= a_n(\mu_{T_1,\ldots,T_{n+1}}),
\end{equation}
and 
\begin{equation}\label{pol10}
\mu_{T_1,\ldots, T_{n-1},\alpha T_nT_{n+1}}= m_n^{(\alpha)}(\mu_{T_1,\ldots,T_{n+1}}).
\end{equation}
\end{lemma}

\proof The proof is based on Lemma~\ref{pol1} and the fact that by \eqref{add3} and \eqref{pol6}, for any Borel set $B\subseteq\C$ we have:
\begin{equation}\label{pol8}
P_{T_n+T_{n+1}}(B)=P_{T_n,T_{n+1}}(1^{-1}(B))
\end{equation}
and 
\[
P_{\alpha T_nT_{n+1}}(B)=P_{\alpha T_n,T_{n+1}}(m^{-1}(B))= P_{T_n,T_{n+1}}((m^{(\alpha)})^{-1}(B)).
\]
In order to prove \eqref{pol7}, consider arbitrary Borel sets $B_1,\ldots,B_n\subseteq \C$. We must show that
\[
\mu_{T_1,\ldots,T_{n-1},T_n+T_{n+1}}(B_1\times\cdots\times B_n)=\mu_{T_1,\ldots,T_{n+1}}(a_n^{-1}(B_1\times\cdots\times B_n)),
\]
i.e. that
\begin{equation}\label{pol9}
\mu_{T_1,\ldots,T_{n-1},T_n+T_{n+1}}(B_1\times\cdots\times B_n)=\mu_{T_1,\ldots,T_{n+1}}(B_1\times\cdots\times B_{n-1}\times a^{-1}(B_n)).
\end{equation}
But by Lemma~\ref{pol1} and by \eqref{pol8},
\begin{eqnarray*}
P_{T_1,\ldots,T_{n-1},T_n+T_{n+1}}(B_1\times\cdots\times B_n)&=& P_{T_1}(B_1)\wedge\cdots\wedge P_{T_{n-1}}(B_{n-1})\wedge P_{T_n+T_{n+1}}(B_n)\\
&=& P_{T_1}(B_1)\wedge\cdots\wedge P_{T_{n-1}}(B_{n-1})\wedge P_{T_n,T_{n+1}}(a^{-1}(B_n))\\
&=& P_{T_1,\ldots,T_{n+1}}(B_1\times\cdots\times B_{n-1}\times a^{-1}(B_n)),
\end{eqnarray*}
and this proves \eqref{pol9}. \eqref{pol10} follows in a similar way. $\endproof$

\vspace{.2cm}

\begin{lemma}\label{permutation}
Let $n\in\N$, and let $\sigma\in S_n$ (the group of permutations of $\{1,2,\ldots,n\}$). Then for any $n$-tuple $(T_1,\ldots,T_n)$ of commuting operators in $\CM$,
\begin{equation}
\mu_{T_{\sigma(1)},\ldots,T_{\sigma(n)}}=\sigma(\mu_{T_1,\ldots,T_n}),
\end{equation}
where identify $\sigma$ with the corresponding permutation of coordinates $\C^n\rightarrow \C^n$.
\end{lemma}

\proof This follows easily from Theorem~\ref{commuting}. $\endproof$

\vspace{.2cm}

\begin{lemma}\label{add one} For $n\in\N$ and $1\leq i \leq n$ define $f_i:\C^n\rightarrow \C^{n+1}$ by
\[
f_i(z_1,\ldots,z_n)= (z_1,\ldots,z_n,z_i).
\]
Then for $n$ commuting operators $T_1,\ldots,T_n\in\CM$, one has that
\begin{equation}\label{add1}
\mu_{T_1,\ldots,T_n,T_i}= f_i(\mu_{T_1,\ldots,T_n}).
\end{equation}
\end{lemma}

\proof Given Borel sets $B_1,\ldots,B_{n+1}\subseteq\C$ we must show that
\begin{equation}\label{add2}
\mu_{T_1,\ldots,T_n,T_i}(B_1\times\cdots\times B_{n+1})= \mu_{T_1,\ldots,T_n}(f_i^{-1}(B_1\times\cdots\times B_{n+1})).
\end{equation}
Clearly,
\[
f_i^{-1}(B_1\times\cdots\times B_{n+1}) = B_1\times\cdots\times (B_i\cap B_{n+1})\times\cdots B_n
\]
so that the right-hand side of \eqref{add2} is 
\[
\begin{split}
\tau(P_{T_1}(B_1)\wedge\cdots\wedge &P_{T_i}(B_i\cap B_{n+1})\wedge \cdots\wedge P_{T_n}(B_n))=\\ &\tau(P_{T_1}(B_1)\wedge\cdots\wedge P_{T_i}(B_i)\wedge P_{T_i}(B_{n+1})\wedge \cdots\wedge P_{T_n}(B_n)).
\end{split}
\]
But this is exactly the left-hand side of \eqref{add2} and we are done. $\endproof$

\vspace{.2cm} 

We will not give the proof of Theorem~\ref{polynomial} in full generality but rather, by way of an example, illustrate how it goes. Consider for instance 3 commuting operators $T_1,T_2,T_3\in\CM$ and the polynomial $q\in\C[z_1,z_2,z_3]$ given by
\begin{equation}
q(z_1,z_2,z_3)= 1+2z_2^2+z_1z_2z_3.
\end{equation}
At first define $\phi_1:\C^3\rightarrow\C^5$ by
\[
\phi_1(z_1,z_2,z_3)=(z_2,z_2,z_1,z_2,z_3).
\]
By repeated use of Lemma~\ref{permutation} and Lemma~\ref{add one} we find that
\[
\mu_{T_2,T_2,T_1,T_2,T_3}=\phi_1(\mu_{T_1,T_2,T_3}).
\]
Next define $\phi_2:\C^5\rightarrow \C^2$ by
\[
\phi_2(z_1,\ldots,z_5)=(2z_1z_2,z_3z_4z_5),
\]
and by repeated use of \eqref{pol10} and Lemma~\ref{permutation} conclude that
\[
\mu_{2T_2^2,T_1T_2T_3}=(\phi_2\circ\phi_1)(\mu_{T_1,T_2,T_3}).
\]
With $\phi_3:\C^2\rightarrow\C$ given by
\[
\phi_3(z_1,z_2)=z_1+z_2
\]
we now have (cf. \eqref{pol7}) that
\[
\mu_{(q-1)(T_1,T_2,T_3)}=(\phi_3\circ\phi_2\circ\phi_1)(\mu_{T_1,T_2,T_3})=(q-1)(\mu_{T_1,T_2,T_3}).
\]
It is now a simple matter to show that for all $\lambda\in\C$,
\[
\tau(\log|q(T_1,T_2,T_3)-\lambda\unit|)=\int_\C \log|z-\lambda|\,\d q(\mu_{T_1,T_2,T_3})(z),
\]
and then by Brown's characterization of $\mu_{q(T_1,T_2,T_3)}$, $\mu_{q(T_1,T_2,T_3)}=q(\mu_{T_1,T_2,T_3})$, as desired.

{\small

\end{document}